\input amstex
\documentstyle{amsppt}
\input epsf
\hoffset0.5in\voffset0.6in
\magnification\magstep1\hsize6truein\vsize8truein
\rightheadtext\nofrills{A. Grassi}
\leftheadtext\nofrills{Divisors and the superpotential}
\topmatter
\title Divisors on elliptic Calabi-Yau 4-folds and the superpotential in 
F-theory, I \endtitle
\author A. Grassi \endauthor
\address   Department of Mathematics, University of Pennsylvania, Philadelphia
PA 19104.
 \endaddress
\NoBlackBoxes
\thanks Research supported in part by NSF grant DMS-94-01495
\endthanks
\endtopmatter
\def\tablerule{\noalign{\hrule\hrule}}

\def\Gr{\operatorname{Gr}}
\def\dim{\operatorname{dim}}
\def\codim{\operatorname{codim}}
\def\Pic{\operatorname{Pic}}
\def\NE{\operatorname{NE}}
\def\bp{\Bbb P^1}
\def\bpt{\Bbb P^1 \times \Bbb P^1}

\def\bda{\bold\Delta}

\def\kda{ K_B + \bold\Delta}

\def \g{\Gamma }

\def\pxb{\pi : X \to B}
\def\C{\Cal}
\def\F{\Cal F}
\def\Cox{\Cal O_X}
\def\Cod{\Cal O_D}
\def\Cob{\Cal O_B}
\def\Coc{\Cal O_C}
\def\d{\Cal D}
\def\c{\chi}
\def\p{\pi ^*}
\def\nxb{\overline {\NE(X/B)}}
\document

\centerline{ABSTRACT} 

\medpagebreak

 Each smooth elliptic Calabi--Yau 4-fold determines both a
three-dimensional physical theory (a compactification of ``M-theory'') and
a four-dimensional physical theory (using the ``F-theory'' construction).
A key issue in both theories is the calculation of the ``superpotential''
of the theory, which by a result of Witten is determined by the divisors
$D$ on the 4-fold satisfying $\chi({\Cal O}_D)=1$.  We propose a systematic
approach to identify these divisors, and derive some criteria to determine
whether a given divisor indeed contributes.  We then apply our techniques
in explicit examples, in particular, when the base $B$ of the elliptic
fibration is a toric variety or a Fano 3-fold.

When $B$ is Fano, we show how divisors contributing to the superpotential
are always ``exceptional" (in some sense) for the Calabi--Yau 4-fold $X$.
This naturally leads to 
certain transitions of $X$, that is birational transformations to
a singular model (where the image of $D$   no longer contributes) as well as 
certain smoothings of the singular model. The singularities which occur
are ``canonical", the same type of singularities of a (singular) Weierstrass 
model. We work out the transitions.
If a smoothing exists, then
the Hodge numbers change. 

 We speculate that divisors contributing to the superpotential are always
``exceptional" (in some sense) for $X$, also  in $M$-theory. 
In fact we show that this is  a consequence of the (log)-minimal model 
algorithm 
in dimension 
 $4$, which is still conjectural in its generality, but it has been worked out 
in various
cases, among which toric varieties.

\bigskip


\smallpagebreak

\it 1991 MSC \rm 14 99, 14 D 99

\newpage

\S 0. Outline of the paper:

\medpagebreak

The general framework for this paper is given by the string  
theories and the
dualities among them (see for example [S], [FL], [V2]).
The original motivation and various applications of this work come  
from physics, while the techniques used
 are in the realm of algebraic geometry.

 We rely in fact on work of Witten, who in [W]  gives necessary and  
sufficient conditions, in mathematical terms,
for the objects of this research, the  divisors  contributing  to  
what is known as
\it the non-perturbative  superpotential in $F$-theory \rm
(see also [BLS, BLS2], [CL], [DGW], [G2], [KLRY], [KS], [KV], [Moh], [My],
[W2]).

 We propose a systematic approach to identify these (smooth)  
irreducible divisors and show how this leads to  questions in  
(birational)
algebraic geometry.

$F$ theory, introduced by Vafa [V], exploits the non-perturbative $SL(2,
\Bbb Z)$ 
symmetry of type $IIB$ string theory in order to produce new types of
physical models associated with elliptic fibrations.  These $F$ theory
models can be regarded as string theories which have been
``compactified" on varieties which admit an elliptic fibration,
often assumed to have a section;
 the $SL(2, \Bbb Z)$ is identified (under the duality between $F$  
and $IIB$ theory) with the symmetry of the homology of the generic fiber.
Our results are phrased in the
context
of $F$-theory (nevertheless, many of the properties stated here are  
also true in a related theory known as $M$-theory).

We thus consider a smooth elliptic Calabi-Yau $4$-fold $\pi: X \to B$ with a 
section, without loss of generality we can assume that $\mu :  X \to W$ is
the resolution
of a Weierstrass model $ \pi _0: W \to B$ (1.1) and that $B$ is uniruled.

 Each smooth elliptic Calabi--Yau 4-fold determines both a
three-dimensional physical theory (a compactification of ``M-theory'') and
a four-dimensional physical theory (using the ``F-theory'' construction).
A key issue in both theories is the calculation of the \it ``superpotential'' 
\rm
of the theory, a sum
$$S(z) = \sum _{D } \exp{<c(D),z>,}$$
 over certain  smooth complex divisors $D\subset X$, where $X$ is a smooth 
Calabi-Yau fourfold  and $z \in H_2(X, \Bbb Z)$; a necessary condition for $D$  
to contribute to the
superpotential is  $\chi(D)=\chi(D, \Cod)=1$ [W].
 
\medpagebreak
In $F$-theory, a divisor $D$ contribute to the superpotential
only if it is ``vertical", that is $\pi(D)$ is a proper subset of $B$;
if the fibration is equidimensional, then such divisors are either components 
of the singular
fibers
(in this case $W$ is necessarily singular), or of the form
$D= \p (C)$, for some smooth divisor $C$ on the threefold $B$ (\S 1).
 We observe that the divisors of the first type are \it ``exceptional" \rm for 
$\mu$
 (in a sense defined precisely in 1.2 and \S 6), are always \it finite \rm in 
number, 
 and can be analyzed by using ``ad hoc" methods,
starting with Kodaira's analysis of the singular fibers and exceptional
divisors of Calabi-Yau $4$-folds.  This is the approach of [KV].
  We study here divisors of the second type. 
\bigskip

  In particular we focus on 2 questions; namely determining when the number of 
such divisors is
finite  and when $D$ is the exceptional divisor of a birational
morphism, which seems to be the case in most examples [W], [BLS], [My], [KLRY],
 [Moh].

 We show how these questions, which are of interest in physics, naturally lead 
to other (open) questions in birational
algebraic geometry. For example, if the log-minimal-model conjecture is true, 
then the divisors contributing to the superpotential are always exceptional, in 
some sense (\S 6, (5.3) and (1.2)).

\smallpagebreak

 We study in detail the case of $B$ Fano (that is $-c_1(B)$ is
very ample): we give an explicit description of all the
divisors contributing to the superpotential (\S 4, 1) and
of the  birational transformations of the Calabi-Yau $4$-folds
which contract these divisors.
 In the Tables (\S 7) we combine these (and other) results.
 What follows is a description of each section:
\bigskip

In \S 1 we
describe properties of such divisors $C$ which 
determine whether $D= \p (C)$ contributes to the superpotential.

\smallpagebreak

In \S 3 we describe our strategy for a systematic approach and develop an 
algorithm.
The fundamental observation [DGW] is that these divisors cannot be \it nef, \rm
 that is there is an effective curve
$\g$ on $B$ such that $C \cdot \g <0$. In particular, it follows from results
of Mori, Koll\'ar, Kawamata that $C \cdot A <0$, where $A$ is the (homology)
class of an effective curve on an ``extremal
ray"  of  the cone of effective curves of $B$ (the dual of the K\"{a}hler 
cone).
 In \S 2 we define the cone
of effective curves (the \it ``Mori cone"), \rm extremal rays, and  properties 
which are relevant in our
 set up.
These objects are, in fact, also the building blocks of the ``Minimal model 
program"
 which, loosely speaking, is an algorithm to construct
a ``preferred" minimal model birationally equivalent to a given
variety. It exactly by following the steps of the algorithm
that we can show that the divisors contributing to the superpotential are, in 
some sense, ``exceptional" (which was hypothesized in [KLRY] and [W]).
We will return to this and the related birational transformations 
in \S 5, 6.
 
 \smallpagebreak

Our strategy consists in examining each extremal ray of the Mori cone and  
argue whether there exists an
effective non-nef smooth divisor $C$ such that $C \cdot \g <0$, with $\g \in 
[R]$. This first step identifies all the possible non-nef divisors. 
Using the technical Lemmas of \S 1 we can then determine the ones with the
right numerical properties to contribute to the superpotential.

This  gives  a straightforward algorithm
which can be applied any time the extremal rays of
the effective cone of $B$ are generated by effective curves; for example, when 
$B$ is Fano, or $B$ is toric, or a $\Bbb P^1$
 bundle over certain surfaces. These cases 
are frequently considered as the basis of Calabi-Yau elliptic fibrations 
[KLRY, Moh, My].
A byproduct of the above algorithm is a list of the fibrations
$B \to S$, with general fiber isomorphic to
$\Bbb P^1$ and  $S$ a surface. This is relevant from the
point of view of the $F$-theory-heterotic duality.

\smallpagebreak

\S 4 contains various  examples.
In particular we concentrate on an equidimensional
elliptic Calabi-Yau $X \to B$, with $B$
a Fano threefold ($c_1(B)>0$) and show that the number of the divisors 
contributing 
to the superpotential is always finite. We use Mori-Mukai's classification of 
such 
threefolds
to describe the divisors of type $D=\p (C)$ which contribute to the 
superpotential 
for each Fano $B$, as well as the $\bp$-fibrations (if any) of $B$ (\S 7).
We also compute the topological Euler characteristic of
$X$, when $X=W$, its smooth Weierstrass model. 

\smallpagebreak

In \S 5 we will show how divisors contributing to the superpotential
naturally are associated to faces of the K\"{a}hler cone of $X$ which lead to
another (necessarily singular) birational model of $X$.

 Some of these divisors are indeed defined by birational contractions
(to the Weierstrass model), see (1.2).
If $X=W \to B$ with $B$ Fano,  any $D = \pi ^*(C)$ contributing to the 
superpotential
 determines a birational transformation $\phi: B \to B'$ with
exceptional divisor $C$ (i.e $codim(\phi(C)) \geq 2$).
We  construct a \it flop \rm of $X$
along $D$, and then contract the image of $D$.
 We obtain an elliptic Calabi-Yau $4$-fold (over $B'$),
 with ``canonical singularities"
 (the same type of singularities as of the Weierstrass model of a Calabi-Yau).
I do not know whether it is possible to build a physical model with these
singularities.
 It is possible in various (\S 5) examples to smooth the singularities
and obtain another Calabi-Yau, where there is \it  no longer a contribution to
the superpotential \rm related to $D$. In this case, the Hodge numbers change.

 In \S 6 we speculate that this is always the case, even in $M$-theory. 
 In fact, a generalized (but still conjectural) version of the minimal model 
algorithm 
implies that given any divisor $D$ contributing to the superpotential
 on a Calabi-Yau $4$-fold $X$ (in $M$ or $F$ theory), 
there is a birational model of the fibration $\rho (X)$
such that $\rho (D) $ does not contribute to the superpotential and $\rho(X)$ 
is singular (at least with \it canonical singularities \rm).
This generalized version (the ``log-minimal model program") has been worked 
out in various
cases, among which are  toric varieties.

\medpagebreak

The divisors contributing to the superpotential thus generate reflections
of the K\"{a}hler cone of $X$ in a ``larger" cone. It would be interesting to 
study
 the Weyl group generated by such reflections, and how this is (if at all) 
related
to the heterotic duals and the change of Hodge numbers of the smoothed variety
(as in [MV], [MS]).

It would also be interesting to study the Calabi-Yau $4$-fold which can be 
``connected" by transitions related to divisors contributing
to the superpotentials:
see also [R] and [ACJM], [AKMS].

\medpagebreak

The core of the paper is in sections 3 and 4 and 6.b: the reader should 
probably 
start
with \S 3 and the general strategy, continue with the examples (\S 4, \S 5, \S 
6.b)
and the Tables (\S 7) and use
the first two sections and \S 6.a as a reference.

\medpagebreak

Finally, in writing this paper I had to give the precedence to some topics
over others. The  parts left out will be investigated in a sequel
(in the not too distant future, I hope).

\medpagebreak

ACKNOWLEDGEMENTS: I have discussed various parts of this project with R. 
Donagi, S. Katz  and in particular D. Morrison and E. Witten.
It is a pleasure to thank them all.
 Thanks also to P. Candelas, M. Larsen, A. Ksir and M. Schneier
for all sorts
of useful suggestions.

\vskip 0.5in

\S 1 Technicalities

\vskip 0.2in

The motivation of this paper comes from describing the divisors contributing to
the superpotential in $F$-theory;  in this context our results are 
complete and more satisfactory at the moment.
Nevertheless, many of the properties stated here apply also to $M$-theory.

We start by considering  a smooth elliptic Calabi-Yau $n$-fold $X$ with a 
section;
that is, $K_X = \Cox$, $h^i(\Cox) =0, \ 0 \leq 1 \leq n-1$ and there is a 
morphism $\pxb$ to a  smooth $n-1$-dimensional variety with
general fiber a smooth elliptic curve. Furthermore there exists a morphism
$s: B \to T \subset X$ which is isomorphic to its image (``the section" of 
$\pi$),
with inverse $\pi _ {|T}$.
  (It is actually enough, for many of the applications considered here,
to consider a ``rational" section, that is, the inverse $s$ of $\pi$ is only 
defined
on an open set in $B$).

 We also assume that the elliptic fiber
degenerates over a non trivial divisor in $B$.
 As a consequence, $h^i(B, \Cob) =0 , \ \forall i>0$; if $\dim(B)=2$,
$B$ is rational; if $\dim (B)=3$, $B$ is uniruled.

\smallpagebreak

\proclaim{Definition 1.0} $\pi _0: W \to B$ is a {\bf Weierstrass 
model} if $W$ can be described by the homogenous equation 
$y^2z=x^3 + Axz^2 + Bz ^3$ in the projective bundle $\Bbb P(\C O \oplus  L^2 
\oplus  L^3)$, with $L$ a line bundle on
$B$ and $A$ and $B$ sections of
$-4L$ and $-6L$ respectively.

 If $L= -K_B$, then $K_W \sim \C O_W$. 
\endproclaim

$W$ is often singular; interesting mathematics and physics
arise from the resolutions of singularities see for example
[MV], [BIKMSV], [KV].
 On the other hand if $-K_B$ is very ample and $h^i(B, \Cob) =0 , \ \forall 
i>0$ then $W$ is a smooth Calabi-Yau manifold.
 Many elliptic Calabi-Yau manifolds can be constructed in this way
(see \S 7).

\smallpagebreak

If  $\pxb$ is an elliptic Calabi-Yau with section,
 we can assume (without loss of generality) that $ \pxb$ is
the resolution
of a Weierstrass model $ \pi _0: W \to B$. In fact:

\proclaim{Lemma 1.1} Let $X \to B$ be a smooth, elliptic Calabi-Yau $n-$fold,
with $B$ smooth. Then

\noindent (1.1.1) $K_X= \p (\kda)$, where $12\bda = 
\sum n_i \Sigma _i, \ n_i \in \Bbb N$.
Here $\Sigma _i$ denotes a component of the locus in $B$ where
 the elliptic curve degenerates;
the summation is taken over all such components.

\noindent (1.1.2)  If the fibration $\pi$ has a section $B \to X$, then there
 exists a Weierstrass model of the fibration and a birational
morphism $\mu$ such that $K_X = \mu ^*(K_W)$ (that is
$W$ has ``canonical" singularities) and the following diagram is
commutative:
$$
\alignat3
& {X} &&  @>\mu>> & \ \ \ {W}\\
& _{\pi } \searrow && & \swarrow _{ \pi _0}\\
& &&\ \ B. &&
\endalignat
$$
(See also \S 6.a.)
\endproclaim

\demo{Proof} (1.1.1) and  existence of the Weierstrass model are due
to Nakayama's Theorem 2.1
 [N]; a discussion of the notation can be found in [MV,I].
 A straightforward argument shows that  the
condition $K_X \sim \Cox$ implies that  $K_X = \mu ^*(K_W)$.
\qed \enddemo

In $F$-theory a divisor $D$ contributes to the superpotential only
if it is ``vertical", that is $\pi(D)$ is not the whole $B$ [W].
In this paper we analyze the divisors $D= \pi ^*(C)$, where $C$ is a smooth 
curve on $B$. 
 The motivation is given by the following:

\proclaim{Observation 1.2} Let $D$ be a smooth divisor in $X$, as above, 
contributing to the superpotential. 

 (1.2.1) If $\pi (D)$ is not a divisor on $B$,  $D$ is necessarily 
exceptional
 for $\mu$.

\vskip 0.1in

\noindent If $\pi : X \to B$ is equidimensional, then either:

 (1.2.2) $ D=  \p (C)$, where $ C \subset B $
 is a smooth irreducible
divisor  such that $\pi _D: D \to C$
is  an elliptic fibration, or

 (1.2.3) $\pi (D) = \Sigma _i$ is a smooth component of the 
ramification divisor and either $D$ is exceptional for $\mu$ or
$W$ is singular along a  subset of $\mu (D)$.
\endproclaim

 \demo{Proof} The existence of a section,
equidimensionality and smoothness of $D$ force $\pi(D)$ to be smooth.
 If $D= \pi ^*(C)$, then a simple computation shows that
the fiber over a general point of $C$ is a smooth elliptic curve
and thus $\pi _D: D \to C$ is an elliptic fibration in the sense
of Kodaira (see also [MV,I]; in fact it is enough to consider the Weierstrass 
model).
\enddemo
 
\proclaim{Observation 1.3} Conversely, if $D= \pi ^*(C)$ is a smooth divisor, 
with $C$ smooth,
then $\pi_D: D \to C$ is an elliptic fibration.
\endproclaim
Note that the divisors of type (1.2.1) and (1.2.3) are always \it finite \rm in
 number and \it exceptional\rm, in some sense. 
As this paper was being written
[KV] appeared, where a particular class of divisors of type (1.2.3) are 
studied. They show in particular that under certain hypothesis,
some $D$ are not exceptional for $\mu$ but contribute to the superpotential 
(as in 1.2.3); this is why we write ``exceptional
in some sense" (see also \S 6.a).

\smallpagebreak

 One application of this work is a criterion to determine under which 
conditions 
the divisors of type (1.2.2) are also finite in number and exceptional.
 If $\pi$ is not \it equidimensional, \rm  then $\pi(D)$ might not be smooth, 
when $D$ is smooth
 (see \S 5.3 for an example). It might be that one should consider,
more generally divisors with mild singularities (see also [G]).
 On the other hand, $\pi$ is indeed equidimensional in many of the examples 
considered in F-theory.

 We should also point out that $\c(D)$ and $h^i(D, \Cod)$  are birational 
invariants
 and Nakayama [Na] showed that there
exists a smooth  birationally equivalent elliptic fibration  which  is
equidimensional over the strict transform of $C$.
We plan to discuss this topic in a continuation of this paper.

\medpagebreak

In the rest of this section we study proprieties of the divisors $D$ of type 
(1.2.2).

Our goal is to reduce the calculation on $B$,  by writing $\c (D)$ as an 
expression on $B$. This is particularly useful when
the geometry of $B$ is well known, for example $B$ is a toric or Fano variety.

(In the following $h^k(V, \C L)=0,$ whenever $k<0$.)
 
\proclaim{Lemma 1.4}  Let $X \to B$ be a smooth, elliptic Calabi-Yau $n$-fold,
with $B$ smooth and let $C \subset B$ and $D \subset X$ be
smooth divisors such that $D = \p (C)$. Then
\smallpagebreak

\noindent $(1.4.1) \quad \ h^m(D, \Cod) = h^m(C, \Coc) + h^{n-1-m}(C, \Coc 
(C)) , \ \forall \ 0 \leq m \leq n-1$ 
 
\noindent $(1.4.2)  \quad \ \phantom{h^m(D, \Cod)} = h^m(C, \Coc) + h^{m-1}(C, 
-\bda _{|C}) , \ \forall \ 0 \leq m
\leq n-1 $

\noindent $(1.4.3)  \quad \ \c (D, \Cod) = \c(C, \Coc) + (-1) ^{n-1} \c(C, 
\Coc (C))$

\noindent $\phantom{(1.4.3)  \quad \ \c (D, \Cod) }=
\c(C, \Coc) + (-1) ^{n-1} \c(C, K_C + \bda_).$

\endproclaim

\demo{Proof}

Note that $\pi_{|D}= ^{\text{def}} \pi_D : D \to C$ is an elliptic fibration 
between
smooth varieties; let us set $\bda _C = \bda _{|C}$. By 1.1,
$12\bda _C$ is a line bundle supported on the ramification locus of
the fibration, which is the complement in $B$ of the locus of the
image of the smooth elliptic curves of the fibration.

 Th 2.1, 7.6, 7.7 in [Ko] apply to $\pi _{D} : D \to C$
and we get the short exact sequences:
$$0 \to H^k(C, (\pi _{D})_ *(K_D)) \to H^k(D, K_D)
\to H ^{k-1}(C, K_C) \to 0, \ 0 \leq k \leq n-2 $$

which give 
$$h^k(D, K_D) = h^k(C, \pi_D *(K_D)) + h^{k-1}(C, K_C) ,
 \ \forall \ 0 \leq k. $$

By the adjunction formula and Lemma 1.1.1 the following equalities hold:
$$\aligned K_D&=( K_X +D )_{|D }= (\Cox + D)_{|D}={\p } (C) _{|D}
={\pi_D} ^* (C_{|C})={\pi_D} ^* (\Coc (C))\\
&= \p (K_B + \bda +C) _{|D} = {\pi_D} ^* ((K_B + C)_{|C} + \bda _C) =
{\pi _D}^* (K_C + \bda _{C}).
\endaligned$$
 
The projection formula [H] now gives  ${\pi_D}_ *(K_D)= K_C + \bda _C = 
C_{|C}. $
Note that $C_{|C}=N_{C/B}$ is the normal
bundle of $C$ in $B$.
 
The statement of the Lemma follows from Serre's duality, applied to $V=C$ 
(resp. $V=D$) and
$L= K_C + \bda_C, K_C$ (resp. $L=K_D$).

(Serre's duality:
 $h^m(V,  L) = h^{r-m}(K_V-L),$ where $L$ a line bundle on a smooth 
$r$-dimensional 
variety $V$.)
\qed
\enddemo

Combining the results in the above Lemma we also get the following Corollary,
which will be used in the explicit computations.

\proclaim {Corollary 1.5 } In the hypothesis of the previous lemma, assume that
$\dim (X)=4$. Then:

$$\aligned (1.5.1) \quad & h^0(D, \Cod)= h^0(C, \Coc) , \\
 &h ^1 (D, \Cod)=  h^1 (C, \Coc) + h^2(C,C_{|C}) = h^1 (C, \Coc) + h^2(B,C) \\
&\phantom{ h ^1 (D, \Cod)}=  h^1 (C, \Coc) + h^2(C, K_C + \bda _C)=h^1 (C, 
\Coc) + h^0(C,-\bda_{|C}) , \\
& h ^2 (D, \Cod)=  h^2 (C, \Coc) + h^1(C,C_{|C}) = h^2 (C, \Coc) + h^1(B,C)\\
&\phantom{  h ^2 (D, \Cod)}=  h^2 (C, \Coc) + h^1(C,K_C + \bda _C) =  h^2 (C, 
\Coc) + h^1(C,-\bda_{|C}) , \\
&h ^3 (D, \Cod)= \phantom{h^2 (C, \Coc) .... }   h^0(C,C_{|C}) = h^0(B,C)-1\\ 
  & \phantom{h ^3 (D, \Cod)}=\phantom{h^0(C,-\bda_{|C})} h^2(C, K_C 
+\bda_{|C})=h^2(C,-\bda_{|C})
\endaligned
$$
  Furthermore,
 $$\c (D, \Cod) = -1/2 (K_C + \bda _{C}) \cdot \bda _{C}= 1/2 K_B
\cdot C^2. \tag 1.5.2$$
\endproclaim

\demo{Proof} When $X$ is a $4$-fold, by the Hirzebruch-Riemann-Roch theorem
for a line bundle $L$
on a smooth surface $C$ we have:
$$\chi (C, L) = \chi (\Coc) + \frac{1}{2}{(L -K_C)} \cdot L$$
and obtain the first equality in (1.5.2) by substituting $L=K_C + \bda _C$.

\smallpagebreak

>From the short exact sequences
  $$ 0 \to \Cob (-C) \to \Cob 
\to \Coc \to 0 \text{ and }
  0 \to \Cob \to \Cob (C)\to \Coc (C) \to 0$$ 
we find:
$$ \chi(D,\Cod )= 2 \chi ( \Cob)- \chi (\Cob (-C)) - 
\chi (\Cob (C)) . \tag 1.5.2' $$
 
On the other hand, the Hirzebruch-Riemann-Roch theorem 
for a line bundle $L$
on a smooth $3$-fold $B$  says that:
$$\chi (B, L) = \chi (\Cob) + \frac{1}{6}L^3 - 
 \frac{1}{4}L^2 \cdot K_B +\frac{1}{12}L \cdot (K_B ^2 + c_2).$$

Substituting  this expression in (1.5.2') for $L=C$ and $L=-C$ respectively,
we obtain the second equality.

 The first set of equalities in (1.5.1) are a direct consequence of (1.4).
 The second set follows from (1.5.2'), since $h^i(B, \Cob)=0, \forall \ i >0$.
 \qed
\enddemo

\medpagebreak
 The following corollary is the first application of the above
machinery; it is obvious, but useful in computation.

\proclaim{Corollary 1.6} 

\noindent (1.6.1) If $\bda _C = \Coc$, then $D$ does not contribute to the 
superpotential.

\noindent (1.6.2) If $\bda _C \neq \Coc$, then $h^1(D, \Cod) = h^1(C, \Coc)$.

\noindent (1.6.3) $h^3(D, \Cod)= 0 \Leftrightarrow h^0(B, C)= 1.$
\endproclaim
\demo{Proof}  (1.6.1) If $\bda _C = \Coc$, then $\c(D, \Cod)=0$, by (1.5.1).

(1.6.2)  Recall that $12 \bda$ is an effective divisor. If $\bda _C \neq 
\Coc$, $h^0(\bda _C) \neq 0$ 
would imply  $h^0( -12 \bda _C) \neq 0$ and
thus  $h^0(12\bda _C) = 0$ (if a divisor and its opposite have non-zero 
sections,
then the divisor is necessarily trivial). \qed
\enddemo

\vskip 0.4in

\S 2. Minimal model theory and the superpotential

\vskip 0.2in

The extremal rays in the sense of Mori  are relevant in this case;
we fix some notation and recall some results of Mori (et al.). Standard 
references
are [CKM], [KMM], [Wil].
  $B$ denotes any smooth complex algebraic variety.
 In sections 3 and 4 we will apply the facts stated here to the case of $B$,
the smooth base threefold of an elliptic Calabi-Yau $4$-fold fibration. (See 
also \S 6.)

\medpagebreak

First some definitions.

 \definition{Definition 2.0} By $\NE(B) \subset \Bbb R^ \ell$ we denote the 
cone generated (over $\Bbb R 
_{\geq 0}$) by the
effective cycles of (complex) dimension 1, mod. numerical equivalence;
and  by $\overline {\NE(B)}$  its closure in the finite
dimensional real vector space $\Bbb R^ \ell$ of all cycles of complex 
dimension 1, 
mod. numerical equivalence (see for example, [CKM]).

Note that $\ell = rk (\Pic (B))$, and in the cases we are considering
here $\ell =b_2(B)$, the second Betti  number of $B$. 

\smallpagebreak

Kleiman's criterion [Kl] says that  $D$ is  \it ample \rm if and only if
$D \cdot \Gamma > 0$ for all $\Gamma  \in \overline {\NE(B)}$;
in particular
$\overline {\NE(B)}$ is
the dual of the  closure of the ample cone
with the duality is given by the intersection pairing between curves and 
divisors.

 The closure of the ample cone is called the \it nef \rm cone:
 a divisor $D$ is  \it nef \rm if and only if
$D \cdot \Gamma \geq 0$ for all the effective curves $\Gamma$ on $B$.
\enddefinition

\smallpagebreak

 The description of $\NE(B)$ for many varieties $B$ can be found in [CKM],
[KMM], [Ma], and for Fano threefolds (the case considered in (4.5)) in [Mt].
We present here some examples that will be relevant in \S 4, in the
description of divisors contributing to the superpotential:

\proclaim{Example 2.1} If $rk (\Pic (B))=1$, then $\NE(B)$ is the positive
 real half-line.$\diamondsuit$
\endproclaim
\proclaim{Example 2.2} If $B= \Bbb F_n$ is a ruled rational surface
$B \to \Bbb P^1$, then $\NE(B)$ is the convex cone generated by $\{f, \sigma _ 
\infty \}$, where $f$ is the (class of the) fiber of the fibration
and by $\sigma _\infty$ the (class of the) unique section with ${\sigma 
_\infty }^2= -n$. $\diamondsuit$\endproclaim

\proclaim{Example 2.3} If $B= \Bbb P^1 \times S$ and $\NE(S)$ is generated by
$\{f_i \}$, then $\NE(B)$ is generated by $\{\ell, f_i \times t \}$, where
$\ell$ is the class of the fiber of the projection $ B \to S$ 
(which is a smooth $\Bbb P^1$) and $t$ is a point in $\Bbb P^1$. 
$\diamondsuit$\endproclaim

\medpagebreak

\medpagebreak

 \definition{Definition 2.4} $R$ is called a \it negative extremal ray \rm on 
the smooth 3fold $B$, 
if $R$ is  an extremal ray of the cone $\overline {\NE(B)}$ in the
usual sense, and 
$K_B \cdot A <0$, for a curve (equivalently, for all
curves) $A$ with homology class spanning  the extremal ray $R$.
 We will write $A \in [R]$.
\enddefinition
\proclaim{Example/Theorem 2.5 (Mori)} If $-K_B$ is ample,
 then every extremal ray is negative and $\overline \NE(B)$ is the convex cone 
generated by
 the extremal rays.
 Furthermore the extremal rays are finite in number.$\diamondsuit$
\endproclaim

 \proclaim{Contraction Theorem (see for example, [CKM], [Wi]) 2.6}  If $R$ is 
a negative extremal ray, then
 there exists a morphism
$\phi _R : B \to B_R$, where $B_R$ is a projective variety
with ``mild" singularities (which can be described) and an irreducible curve 
$E \subset B$ is contracted
by $\phi _R$ if and only if the homology class of $E$
belongs to the extremal ray $R$. Furthermore $rk (\Pic (B)) > rk (\Pic (B_R))$;
 $\phi _R$ is called the contraction morphism.
\endproclaim

(The singularities which occur are called \it ``terminal" \rm ; if $\dim 
(B)=2$ these are the smooth points.)
\remark{Remark 2.7}
In general, if any morphism contracts a curve on one extremal ray, then it
necessarily contracts all the effective curves on the same 
extremal ray.
\endremark

\example{Examples (Contraction morphisms and extremal rays) 2.8 }  In example 
(2.1), if $K_B \cdot \g <0$, for an effective curve $\g$ on $B$, 
(equivalently, 
all effective curves)  then $\NE(B)$ consists of one negative extremal ray 
and
the corresponding contraction morphism sends $B$ to a point.
If $K_B \cdot \g \geq 0$, then there is no negative extremal ray.

In example (2.2), $f$ is always a negative extremal ray ($K_B \cdot f=-2$) and 
the corresponding contraction 
morphism  gives the structure of $\Bbb P^1$-bundles
$\Bbb F_n \to \Bbb P^1$. 

 On the other hand,  $K_B \cdot \sigma _\infty  <0$ only when $n=1$; in this 
case $\sigma _\infty$ is the only negative extremal
ray.
 The corresponding contraction morphism is $\Bbb F_n \to \Bbb P^2$ the blow up
of $\Bbb P^2$ at a point.
Note that we can always contract $\sigma _\infty$, independently of the value 
of 
$n$.
The image surface however will always be singular unless $n=1$.
In fact the ``mild" singularities mentioned above (in the statement
of the Contraction Theorem)
are exactly the smooth points
when $\dim(B)=2$.
$\diamondsuit$
\endexample

 \proclaim{Example/Theorem 2.9 (Mori)} If dim $B=3$, then the exceptional 
locus $C_R$ of a birational morphism $\phi _R: B \to B_R$ associated to a 
negative
 extremal ray 
$R$
is one of the following reduced divisors:

\noindent (2.9.1) $C_R$ is a $\Bbb P^1$-bundle over the smooth curve
$\phi _R (C_R)$, 
with

$\frac{1}{2} K_B \cdot {C_R}^2 = 1- g(\phi _R (C_R))$;

\noindent (2.9.2) $C_R \sim 
\Bbb P^2$ with $\C O_{C_R}(C_R) =\C O_{\Bbb 
P^2}(-1)$

\noindent (2.9.3) $C_R \sim \Bbb P^2$ with $\C O_{C_R}({C_R}) =\C 
O_{\Bbb 
P^2}(-2);$

\noindent (2.9.4) ${C_R} \sim \Bbb P^1 \times \Bbb P^1$ with
$\C O_{C_R}({C_R}) =\C O_{\Bbb P^1 \times \Bbb P^1}(-1,-1);$

\noindent (2.9.5) $C_R$ is a singular quadric surface in $\Bbb P^3$.
 $\diamondsuit$

In cases (2.9.1)-(2.9.2) $B_R$ is a non singular $3$-fold;
$\phi_R(C_R)$ is a quadruple point in case (2.9.3) and a double point
otherwise.
\endproclaim
\demo{Proof} See Mori [Mo], [CKM] or [Wi]. \qed \enddemo

\newpage

\S 3. The algorithm.

\vskip 0.2in

We now consider the case of $\pi: X \to B$, an elliptic fibration
of a Calabi-Yau $4$-fold $X$, with $C$ and $D= \p (C)$ smooth divisors.
This happens when $\pi$ is an equidimensional elliptic 
fibration with
section, as we saw in (1.2).

\smallpagebreak

The following remarks are the building blocks of our strategy:

\proclaim{Remark 3.1} If
$D = \p (C)$ contributes to the superpotential, then there exists
an extremal ray  $R$ on $\overline {\NE(B)}$ such that
$C \cdot A <0$, for all the curves
$A$  on the ray $R$. In particular the set of divisors contributing to the 
superpotential is contained in the set of extremal rays.
\endproclaim

In fact, $D$ and $C$ cannot be nef divisors [DGW]; 
$C$ is non-nef 
if and only if $C \cdot A <0$, for all 
$A$  on some extremal
ray $R$ of $\overline {\NE(B)}$. 

\smallpagebreak

\proclaim{Remark 3.2} If $C$ is not nef, all the curves on the extremal ray $R$
 must be
contained in $C$. In particular, {\bf if} there exists a morphism $B \to S$
contracting exactly
 the curves on the extremal ray $R$, then $\dim (B)= \dim (S)$.
\endproclaim

 \proclaim{Strategy} ($\bullet$) We consider cases where the extremal rays of
 $\overline {\NE (B)}$ are generated by effective curves.

($\bullet \bullet$) For  each extremal ray $R$, we determine whether there 
exists an
effective smooth divisor $C$ such that $C \cdot \g <0$, for $\g$ on the ray 
$R$. 

($\bullet \bullet \bullet$) If such a $C$ exists, we check its numerical 
properties.
\endproclaim
 
In most relevant
 cases  (in $F$-theory) this strategy gives a quick algorithm to determine the 
divisors of this form contributing to the superpotential. We will do so 
explicitly in \S 4.

In fact, 
 the extremal rays generate the cone of effective curve
when $B$ is Fano (2.4),  toric [B], [O], or a $\Bbb P^1$
 bundle over certain surfaces. These 
cases are frequently considered as the basis of Calabi-Yau elliptic fibrations 
($\bullet$).

 Often the extremal rays are defined in terms of morphisms 
(see (2.7)); this is in fact always the case  for  negative extremal
rays (by the contraction theorem) ($\bullet \bullet$). 

  At the same time (by looking at the extremal rays) we can also 
describe the $\Bbb P^1$ bundle structure (if any) of $B$. This is relevant 
from the point of view of duality with heterotic theory.
 We will do so explicitly in \S 4 and in the Table.

 We will use (1.5) and (1.6) for ($\bullet \bullet \bullet$).

\smallpagebreak

 If the fibration
$\pi$ is equidimensional, and the K\"{a}hler cone of
$X$ (and hence the cone of $B$) is polyhedral, that is has a finite number of 
extremal rays,
then the number of divisors contributing to the superpotential is \it finite 
\rm (1.2).

\medpagebreak

Another advantage of this approach is that in our examples
we get a map

\smallpagebreak

 $$\{ \text{divisors contributing to the superpotential} \}   \rightarrow \{ 
\text{faces of the K\"{a}hler cone of } X \}.$$ 

\bigpagebreak

Note that the divisors of type (1.2.1) and (1.2.3) are always
associated to a face of
the K\"{a}hler cone of $X$ as they come from resolving the Weierstrass
model of $X$.

 For the divisors of type (1.2.2) the question is more subtle:

\proclaim{Remark 3.3} In terms of the dual ``nef" cone,  the morphism 
associated to a chosen negative extremal
ray gives a divisor class on the boundary of the ``nef cone" of $B$ (2.6).
 When we start from a divisor contributing to the superpotential, this
``face" of the nef cone must lead to another (birational) model of $B$ by 
(3.2).
\endproclaim

\smallpagebreak

 In \S 5 we will show how divisors contributing to the superpotential
 are associated to faces of the K\"{a}hler cone of $X$ which lead to
another (necessarily singular) birational model (of $X$).
 
 We speculate that this is always the case, even in $M$-theory and show how 
this is related to various conjectures in algebraic geometry [M], [K], [KMM].

\medpagebreak

{\bf The case of negative extremal rays}

\medpagebreak

Only the negative extremal rays which determine birational contractions (3.2)
 are relevant for our purposes.

 In this this case ($\dim (B) \leq 3$) there is a unique non-nef divisor $C_R$ 
such that
$C_R \cdot \g <0, \ \forall \ \g \in [R]$  (2.9); we also have a complete
list of the possible $C_R$ which occur.
 We only consider here smooth divisors $C_R$ (see [W]);
the case (2.9.5), the quadric cone in $\Bbb P^3$ should be also of interest,
as it is a divisor of simple normal crossings [G]. This will be investigated 
in a forthcoming paper.
\smallpagebreak

 The following proposition follows directly from
 (1.5) together with Mori's description given in (2.9).

\vskip 0.2in

\proclaim{Proposition 3.4}  Let $R$ be a negative extremal ray associated to a 
birational morphism $\phi _R$, $C_R$  the unique
exceptional divisor  and
$D_R= \p (C_R)$, as in (2.9).
 In cases (2.9.2), (2.9.3), (2.9.4),
 $$ h^0(D_R, \C O_{D_R})=1, \ h^1(D_R, \C O_{D_R})=
h^2(D_R, \C O_{D_R})= h^3(D_R, \C O_{D_R} )=0,$$
and $D_R$ always contributes to the superpotential.

In case (2.9.1),  
$$\chi (D_R)=1/2 K_B \cdot {C_R}^2 = 1- g(\phi _R(C_R)) =1 $$
if and only if $\phi _R({C_R})$ is a rational curve; furthermore
$$\align h^0(D_R, \C O_{D_R})&=1, \ h^1(D_R, \C O_{D_R})= h^1(C, \Coc)=0\\
h^2(C, \C O_{D_R})&= \c( D_R, \C O_{D_R}) -1 -h^1(C, \Coc) = h^3(D_R, \C 
O_{D_R} )=0
\endalign
$$
and $D_R$ contributes to the superpotential if and only if
 $C_R$ is rationally ruled.
\endproclaim

We speculate that these are exactly the cases when $C_R$ deforms
whenever $B_R$ deforms.

\vskip 0.5in

\S 4 Examples (The algorithm at work.)

\medpagebreak

We apply the algorithm outlined in \S 3 to various examples.
Because we restrict ourselves to divisors of type (1.2.2) (that is of the form
$D= \p (C)$, with $D$ and $C$ both smooth), we describe 
on $B$ the relevant divisors $C$ (such
that $D =\p (C)$ contributes to the superpotential).  
If $X=W$, then
these are all the divisors contributing to the superpotential; 
and we can write the superpotential (\S 7).
 The divisors not of this form are always finite in number
(in $F$-theory), are ``exceptional" in some sense (see (1.2)), and can
 be described with other methods.

\smallpagebreak

\proclaim{Example 4.1} If  $b_2(B)=rk(\Pic(B))=1$, no divisor 
of the form $D= \p (C)$ contributes to the
superpotential and there is no fibration $B \to S, \ \dim (S)\neq 0$.
\endproclaim
 In fact, in this case, $\NE(B)$ is a half-line (see 2.1): any divisor
containing all of the line would also contain all of $B$.

 In particular, $-K_B= c_1(B)$ is necessarily an ample divisor:
in fact $-12K_B$ is the effective, non trivial, divisor image (under $\pi$) of 
the singular fibers.
$B$ is a Fano threefold; such varieties were classified by Iskovskih [I1, I2].
 Among those are 
$B= \Bbb P^3$ and $B=Q$ the smooth quadric in $\Bbb P^4$. The complete
list will appear in \S 7. $\diamondsuit$

\smallpagebreak

\proclaim{Example 4.2 (no. 27, Table 3)} If  $B = \bp \times \bp \times \bp$ 
no divisor
of the form $D=\p(C)$ contributes to the superpotential. 
\endproclaim
In fact, $\overline {\NE(B)}$ is a cone with $3$ edges in $\Bbb R^3$: each edge
being a fiber of the projection to $2$ of the factors (see (2.3) and (3.2)).
$\diamondsuit$

\smallpagebreak

\proclaim{Example 4.3 (for $n=1$, no. 28, Table 3)} If $B = \bp \times \Bbb 
F_n$, $n \geq 1$, then no 
divisor
contributes to the superpotential when $n \neq 2$, and 1 divisor
contributes when $n=1$. In the latter case, the divisor is determined by
a negative extremal ray of type (2.9.1).

There is a $\bp$ fibration $B \to \Bbb F_n$ and a $\bp$ fibration $B 
\to \bp \times \bp.$ 
\endproclaim
(This is analysed in [W], for $n=1$.)

In fact, $\NE(B)$ is generated by $\{ \ell, f, \sigma \}$ where
$\ell$ is a fiber of  $p: B \to \Bbb F_n$, $f$ a fiber
of $B \to \bp \times \bp$, and 
$\sigma= \sigma_\infty \times \{t \}, \ t \in \bp$, as in (2.2)-(2.3).

 There is no non-nef divisor associated to 
$f$ or $\ell$ (3.1), (2.8); if we set $C=p^{-1}(\sigma _ \infty)$,
then $C \cdot \sigma = -n$ is non-nef when $n >0$.
 Note that $C \sim \bp \times \bp$ and
$h^0(C, \Coc)=0, \ h^1(C, \Coc)=h^2(C, \Coc)=0$; we identify 
the Picard group of $\bp \times \bp$ with $(a,b)$: a curve is effective
when $a, b \geq 0$.

 We need to compute $h^i (C, C_{|C})$ (for $n >0$) and determine whether $C$
 contributes  to the superpotential, by 1.5.
 A simple computation gives $C_{|C}=(-n,0)$ and immediately
$h^0(C, C_{|C})=0, \ h^2(C, C_{|C})= h^0(C, (n-2,-2))=0$ (by Serre's duality).
 It follows also that $h^1(C, C_{|C}) >0$ if $n \geq 2$ and
$h^1(C, C_{|C}) = 0$, for $n =1$ (Kunneth's  formula).
$\diamondsuit$

\smallpagebreak

\proclaim{Example 4.4 [DGW]} If  $B= S \times \Bbb P^1$,
where $S$ is a general rational elliptic surface with section, then
there are infinitely many divisors contributing to the superpotential,
corresponding to the negative extremal rays of $NE(S)$.
There is one $\bp$ fibration to $\bp \times \bp$ and one to $S$.
\endproclaim

  It is not hard to see that the generators for $\overline {\NE(S)}$ are 
$f$ and $s_\alpha$, where $f$ is a fiber and $s_\alpha$ a section of the
elliptic fibration $p: S \to \Bbb P^1$.
 We can choose $S$ so that the $s_\alpha$'s are infinitely many;
it turns out that every $s_\alpha$ is a a negative extremal ray.

 Then $NE(B)$ (see (2.3)) is generated by $\{f \times t, s_\alpha \times t, 
\ell \}$, with $s _\alpha, f$ as above, 
$t \in \Bbb P^1$, and $\ell$ a fiber of $B \to S$.

  $f \times t$ and $\ell$ do not determine  divisors contributing to the 
superpotential: they are in fact  the general fibers of
$B \to \Bbb P^1 \times \Bbb P^1$ and $B \to S$ respectively
(3.2), while $s _\alpha \times t $ is a negative extremal ray  on
$B$, for all $\alpha$. 
The corresponding divisor $C_{s_\alpha}$ is of type (2.9.1), 
 is isomorphic to $ \Bbb P^1 \times \Bbb P^1$ and contribute to the 
superpotential. It follows that the divisors contributing to the 
superpotential are exactly the 
$C_{s_\alpha}$ (3.1).
$\diamondsuit$

\smallpagebreak

\proclaim{Example/Theorem 4.5.1} Assume that $B$ is a Fano variety, that is
$c_1(B)= -K_B$ is ample. Then 
the divisors $D= \p(C)$ contributing
to the superpotential are exactly the  
exceptional divisors $C=C_R$ (corresponding to the contraction
of extremal ray $R$) listed below:

\noindent (4.5.1.1) $C_R$ is a $\Bbb P^1$ bundle over a smooth rational curve 
($\phi _R (C_R)$),

\noindent (4.5.1.2)-(4.5.1.3) $C_R \sim 
\Bbb P^2$,

\noindent (4.5.1.4) ${C_R} \sim \Bbb P^1 \times \Bbb P^1$.
\endproclaim

\smallpagebreak

\demo{Proof} If $B$ is Fano,  the negative
extremal rays generate $\overline {\NE(B)}$;
thus $C$ is non-nef only if it contains all the
curves $E$ in the homology classes of some 
negative extremal ray $R$.
Mori's result says that all such curves
span the exceptional
locus of the morphism $\phi_R$, described in 2.9.
 Then necessarily,
$C=C_R$  and $\phi_R$ is a divisorial
contraction.  The statement follows from (3.4). \qed \enddemo

\proclaim{Corollary 4.5.2} Let $B$ be a Fano $3$-fold.

\noindent (4.5.2.1) If $\pi$ is equidimensional, then
there is  only a finite number of divisors contributing
to the superpotential. 

\noindent (4.5.2.2) If $X=W$ (the Weierstrass model is smooth), then 
 the divisors contributing to superpotential
are exactly the exceptional divisors of Mori contractions listed above.
\endproclaim

\demo{Proof}
If $B$ is Fano (by the``Cone theorem") there are only finitely many
negative extremal rays, hence a finite
number of such divisors $D$ on $X$ contributing
to the superpotential. \qed \enddemo

In the examples in [W] Witten shows that ``a superpotential is not generated by
instantons by showing that any divisor $D$ on $X$ has
$\chi (D) \neq 1$, or we show that a superpotential is generated by showing
that some choice of the cohomology class there is precisely one complex 
divisor $D$,
which moreover has $h_1=h_2=h_3=0$".

This is exactly what always happens for $B$ Fano:

\proclaim{Corollary 4.5.3}  Let $B$ be a Fano $3$-fold and $D= \pi ^*(C)$.
Then either $\chi (D) \neq 1$ or $h^0(D)=1, \ h^i(D)=0, \ i \neq 0$.
\endproclaim

 Mori-Mukai [MM] classified all Fano  threefolds and Matsuki [Mt]
described the extremal rays for each of them: 
 the relevant divisors are the one corresponding to birational contractions 
(3.2).
We apply our algorithm to each case in their list and
we determine the divisors of type (1.2.2) contributing the the superpotential 
(\S 7).
 The only delicate point is when $C$ is of type (2.9.1), i.e.
a $\bp$-bundle over a smooth curve $L = \phi _R (C)$: $C$ contributes if and 
only if $L$ is rational (3.4).
 The following identity is useful to compute $g(L)$ the genus of $L$ (notation
as in (2.9):

$${(-K_B)}^3={(-K_{B_R})}^3 -2 \{ -K_{B_R} \cdot L - g(L) +1 \}.$$

\definition{Example 4.5.4} ($b_2=3$, no.9 in [MM].)
$B$ is the blow up of the cone over the Veronese surface  $ R_4 \subset
\Bbb P^5$
with center a disjoint union of the vertex and the quartic in
$R_4 \sim \Bbb P^2$.
(Recall that the Veronese surface is $\Bbb P^2$ embedded in $\Bbb P^5$
by its linear system of conics.)

The Matsuki-Mori-Mukai classification says that $\overline {\NE(B)}$
is generated by $4$ curves (the extremal rays):
$R_1$, the ruling of the exceptional divisor over the quartic,
$R_2$, the strict transform of a line in the Veronese surface,
$R_3$, the ruling of the exceptional divisor which is the strict transform
of the ruling over the quartic,
$R_4$, a line in the exceptional divisor of the blow up of the vertex of 
the Veronese cone.
  Furthermore the corresponding extremal contractions $\phi _{R_i}: B \to
B_{R_i}$
are all birational:
$R_1$ and $R_3$ are of type (2.9.1), while $R_2$ and $R_4$ are of type (2.9.3).
The exceptional divisors of $\phi _{R_i}, \ i =2, 4$ contribute
to the superpotential, while the others  are $\bp$-bundles
over a curve of positive genus (the plane quartic) and do not contribute
(4.5.1).

 $B_{R_i}, \ i =1, 3$ is the blow up of the Veronese cone
with center a plane quartic; while  $B_{R_i}, \ i =2 ,4$ is isomorphic
to the blow up of the Veronese cone with center the vertex ($b_2=2$, no. 36). 
The extremal transition with exceptional divisors
contributing to the superpotential lead in this case to a \it
singular \rm variety $B_R$.
\qed
\enddefinition

\vskip 0.4in

\S 5. Transitions of CY $4$-folds I: \it The case of negative extremal rays. 
\rm

\vskip 0.2in

 One of the examples studied by Witten [W] is  $B$ the blown up of $\Bbb P^3$
at a  point (this is no.35 in Mori-Mukai's list). A puzzle arises here:
 while there is no divisor on $\Bbb P^3$ contributing to the superpotential,
the exceptional divisor of the blow up contributes on $B$.

 We will show that all the divisors contributing to the superpotential
are always ``exceptional" in some sense, at least when $B$ is Fano.
The general statement depends on the (log)-minimal model conjecture
which will be discussed in \S 6.

Note that some of these divisors are actually ``defined" by birational 
contractions
(to the Weierstrass model), see (1.2). The divisors considered in 
[KV] are of this form.

\medpagebreak

What will follow is in fact a $4$-dimensional analogue of the construction
in [MV, II]: in that case $B = \Bbb F_1$, $C_R$ is the curve with
self-intersection $-1$ (this curve is in fact a negative
extremal ray, (2.2)). Morrison and Vafa perform a toric flop of the 
holomorphic image of $C_R$ in $X$ ($X \to B$ has  a section), and then contract 
the
image of  $D = \pi ^{-1} (C)$ (which is a Del Pezzo surface). Finally they 
smooth the singularity.

\medpagebreak

Similarly, we consider $\pi: X \to B$, with $X$ equal to its smooth 
Weierstrass 
model, $B$ Fano and assume that $D=\pi^*(C)$ is a divisor contributing to the 
superpotential. Then $C=C_R$ is the exceptional divisor
of the contraction morphism $\phi _R : B \to B_R$ associated to the negative 
extremal ray $R$ (\S 3, (4.5)). $D$ cannot be contracted immediately 
(see 6.8), so (as in [MV, II]) we first start with a ``flop" (5.1)
and follow with a contraction (5.2).

As in [MV] we assume, for simplicity's sake, that
the fibration $\pi$ is general, i.e. that there is
only 1 section.

\medpagebreak

The new threefold $\bar X_R \to B_R$ is elliptically fibred, but it is 
singular: in the cases where $X$ is a general smooth Weierstrass model,
(as in [MV, II]) the singularities can be described precisely.
  These type of singularities are called ``canonical" (\S 6) and are the same 
type of singularities that occur on the \it singular \rm  Weierstrass
models.
It is not clear to me whether  a physical model can be built
with these singularities.

 If the singularity can be smoothed (we explicitly do so in various cases), 
then the resulting Calabi-Yau
will have different Hodge numbers.

 It is an interesting question to investigate this change and how it might
be related to the exceptional divisor contracted (as in [MV,II]).

\smallpagebreak

\proclaim{Proposition (the flop) 5.1} Let $C= C_R$ as in (4.5.1) and
(4.5.2). Then there exists
a contraction
$B \to B_R$, where $B_R$ is another threefold and a birational
transformation (``flop")  $X \dasharrow X_R$ such that the following diagram 
is commutative
$$
\CD
X      @>\mu _R>>     X_R \\
@V{\pi}VV           @VV{ \pi _R }V   \\
B        @>\phi _R>>   B_R \\ 
\endCD 
$$
 $X_R$ is smooth only
if $C_R$ is of type (2.9.1).
\endproclaim

\demo{Proof} (Following Matsuki.)
We have assumed the existence of a section
 of the elliptic fibration; so there exists a smooth threefold $T$ isomorphic 
to
$B$ in $X$ (a ``copy of $B$ in $X$); by $C_R$ we 
will
 denote both the surface in $B$ and its isomorphic image in $T$. 
   We can ``duplicate" the contraction of  $C_R$ in $B$ (2.9) in its
holomorphic image in $T$ (and $X$) and obtain a birational
transformation $X \to Z$  6.8). 
 Matsuki in [Mt] considers a similar situation and  explicitly 
constructs the flop of each surface $C_R$ in $X$, for each
$R$.  The pictures are fairly self-explanatory: 

\noindent $\bullet$ the large ovals
denote $T$, the image of $B$ in $X$ and its images after the 
the blow ups and blow downs,

\noindent $\bullet$ the object in the ovals denote the image of $C_R$ in $X$ 
and their images after the various
birational transformations,

\noindent $\bullet$ the ``parachute type" objects in the  $X$ and $X_R$ denote 
$D$ and its image $D_R$
after the ``flop".

 It is clear from the picture
that $D_R$ has intersection positive with the
fiber of the contraction with $Z$, while $D \cdot R <0$
($R$ is the fiber of the contraction $X \to Z$).
 We have performed a \it ``log-flip" \rm  with respect to
$D$ (see also 6.8).

For a detailed description see [Mt], pages 30-36 and also \S 6.
 
\newpage

\noindent \it $C_R$ is of type (2.9.1): \rm
$C_R$ is $\bp$-bundle over $\phi _R(C)$, and $B_R$ is a smooth
$3$-fold.
   The shaded area is a vertical ``section" of $\mu(D)= D_R$, which is 
isomorphic to the DelPezzo surface
which is obtained by blowing up $\Bbb P^2$ at $8$ points
(see also 5.2, 6.8 and [MV,II]).

\medpagebreak

\noindent\phantom{it'salongway}\epsfxsize=3in\epsfbox{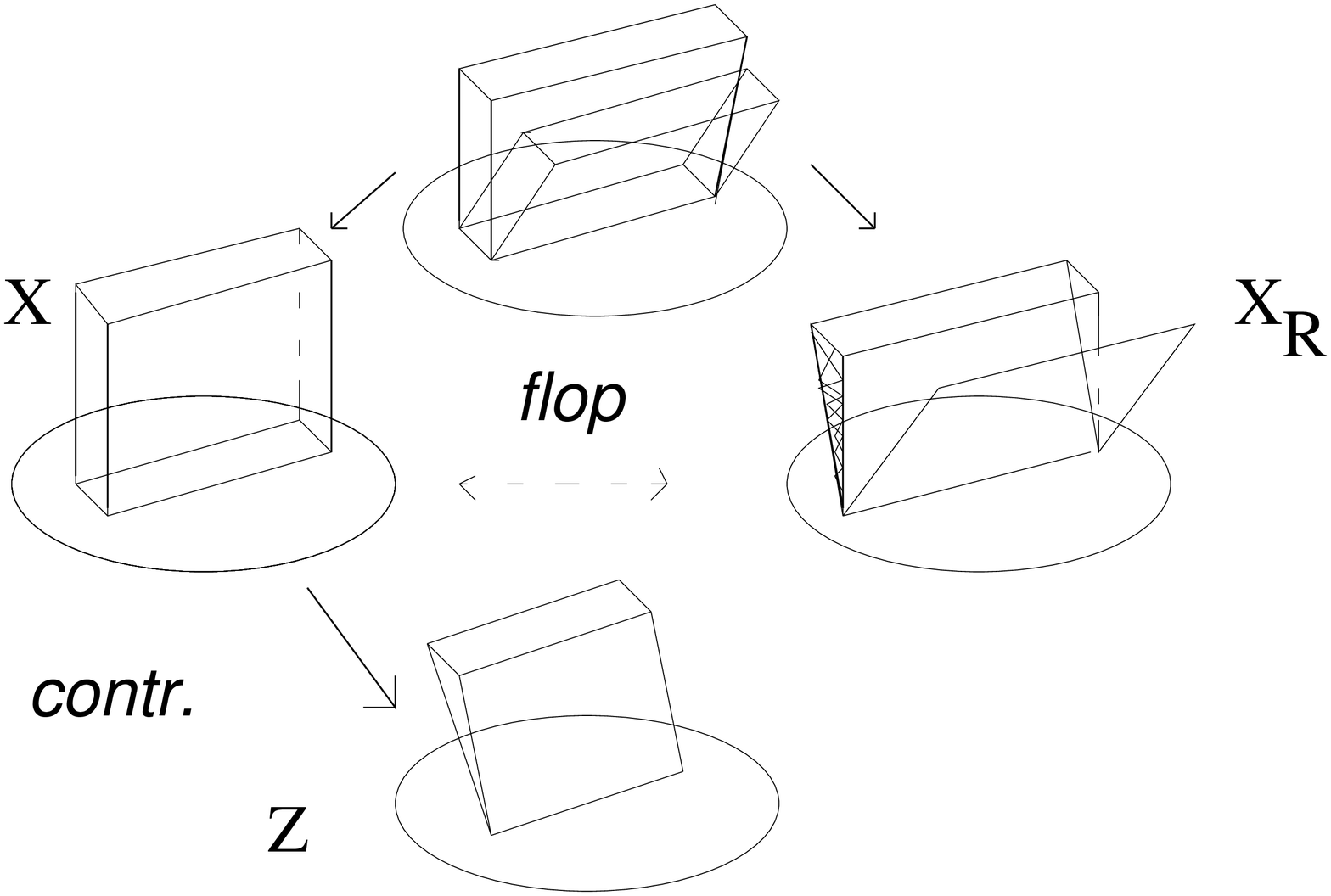}
 
\medpagebreak

$X_R$ is smooth.

\vskip 0.3in

\noindent\it $C_R$ is of type (2.9.2): \rm $C_R \sim \bpt$,
$\phi _R(C_R)$ is an ordinary double point in $B_R$.

\medpagebreak

\noindent\phantom{it'salongway}\epsfxsize=3in\epsfbox{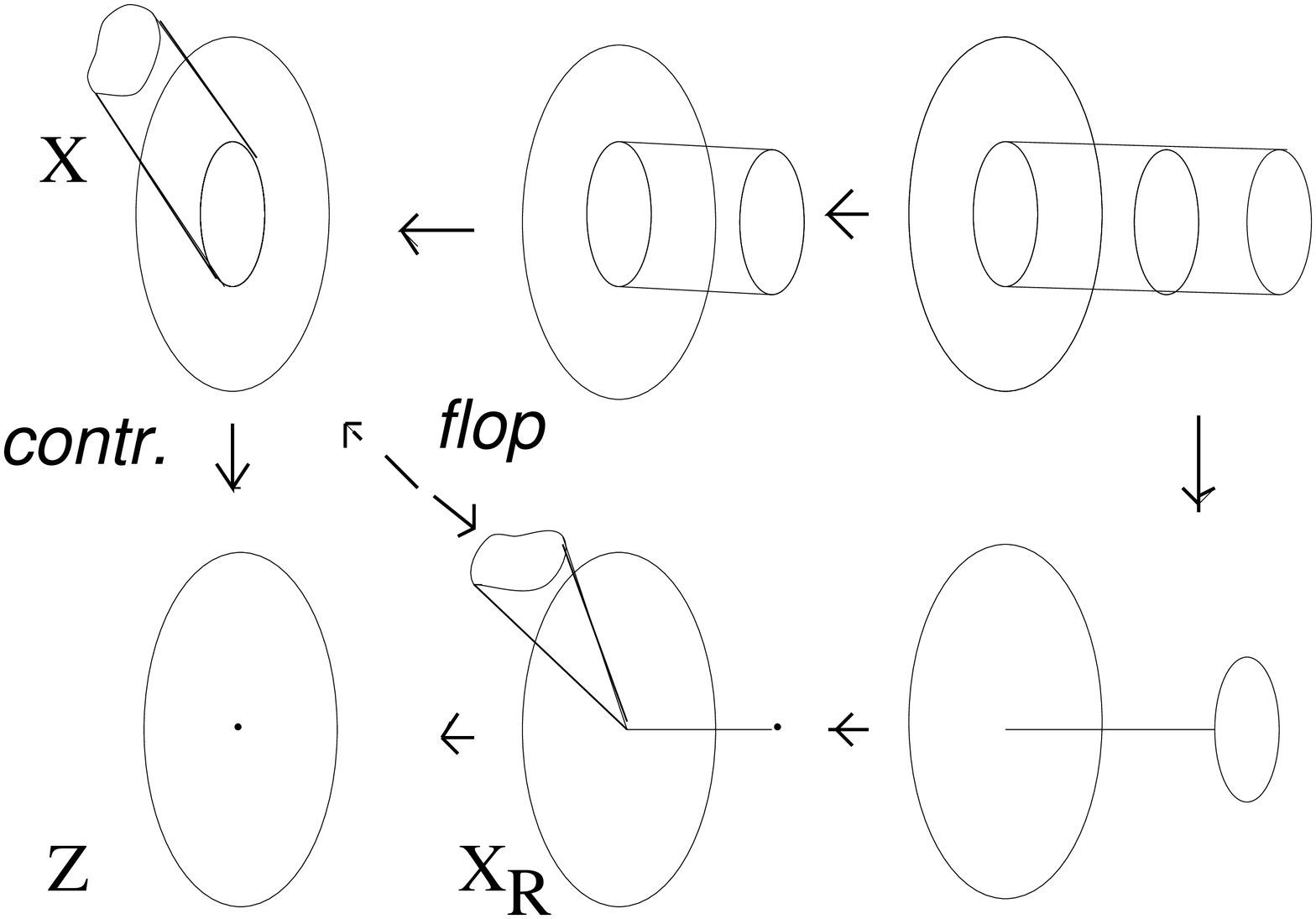}
 
\medpagebreak

 $X_R$ is singular along a $\bp $ (the ``fat" point in the picture).
\vskip 0.3in

\noindent \it $C_R$ is of type (2.9.3): \rm
$C_R \sim \Bbb P^2 $,
$\phi _R(C_R)$ is a quadruple  point in $B_R$.

\medpagebreak

\noindent\phantom{it'salongway}\epsfxsize=3in\epsfbox{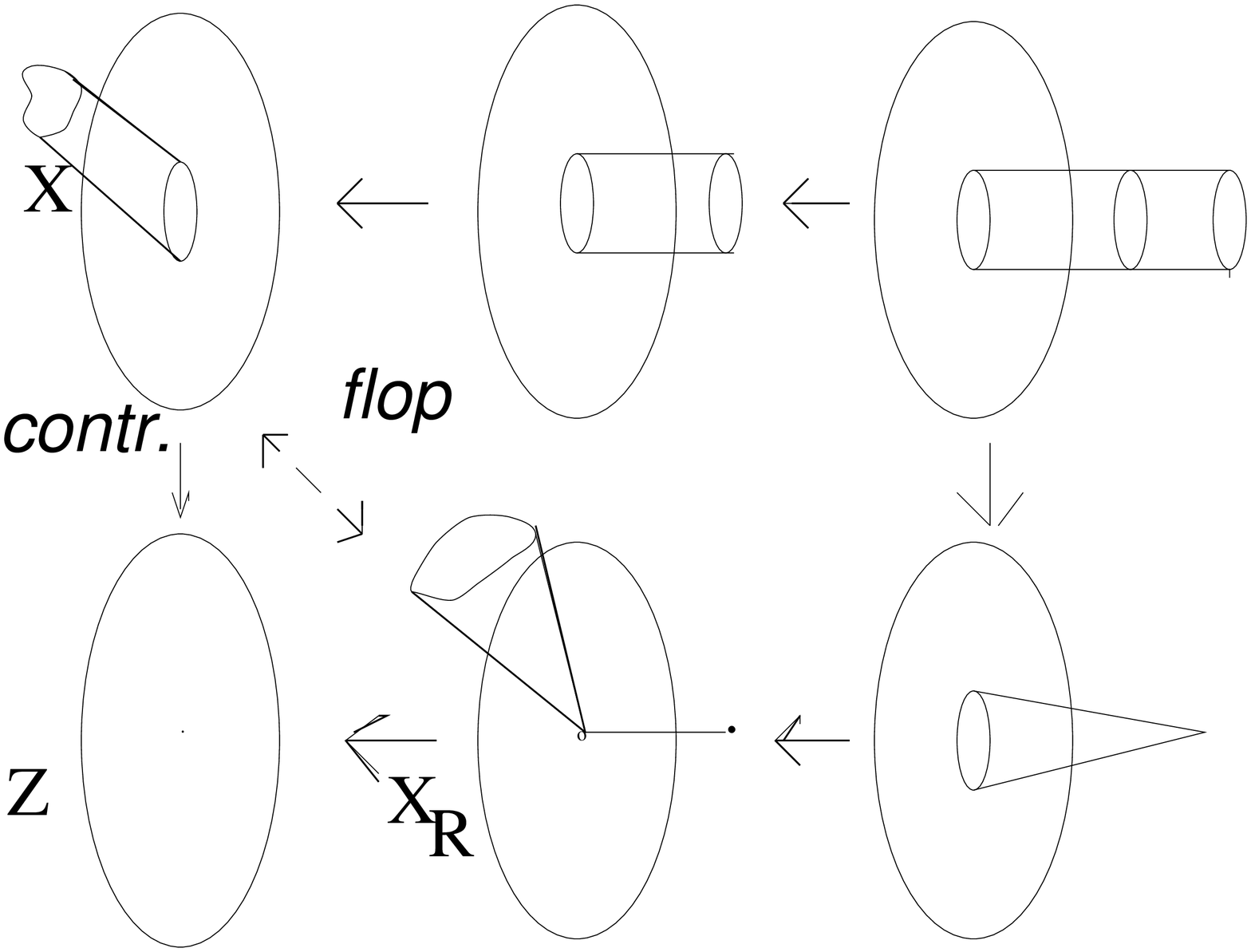}
 
\medpagebreak

 $X_R$ has singular  point (the ``fat" point in the picture).
\vskip 0.3in

\noindent\it $C_R$ is of type (2.9.4): \rm
$C_R \sim \Bbb P^2 $,
 $B_R$ is non-singular.
\medpagebreak

\noindent\phantom{it'salongway}\epsfxsize=3in\epsfbox{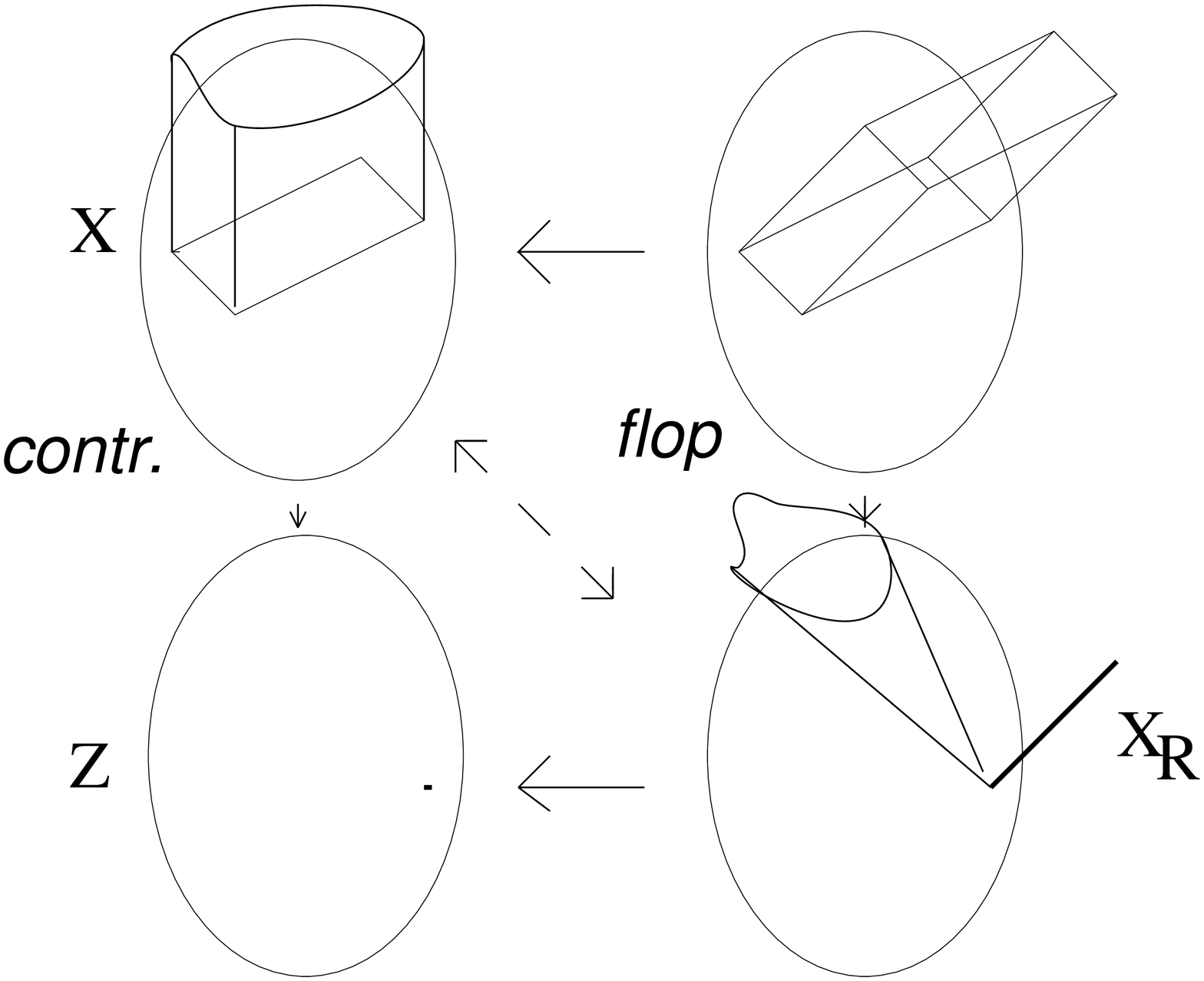}

 $X_R$ has singular  point (the ``fat" line in the picture).
\qed
 \enddemo

As Matsuki points
out these are not the only flops which can occur; however our goal here 
is to show 
that
we can ultimately contract the
image of $D$, which cannot be contracted in $X$
(see \S 6). However it is possible that one would need to consider
 other type of flops to describe all the reflections (and corresponding Weyl
 group)
of the K\"{a}hler cone of $X$, in the enlarged K\"{a}hler cone, determined
by the divisor contributing the superpotential.

Note the flops used above are toric, even when $B$ is not toric.

\proclaim{Proposition (the contraction) 5.2} There is  birational 
transformation  
$\rho: X \to \bar X_R$ with exceptional divisor 
$\rho (D)$ and an elliptic fibration $\bar X_R  \to B_R$ (with section)
such that the following diagram is commutative:
$$
\CD
X      @>\rho>>     {\bar X_R} \\
@V{\pi}VV           @VV{ \bar \pi _R }V   \\
B        @>\phi _R>>   B_R. \\ 
\endCD 
$$

$\bar X_R$ has {\bf canonical} singularities ($\rho ^*(K_{\bar X_R}) = K_{ 
X})$.
\endproclaim

If the Weierstrass model of $X$ is singular, it  has canonical
singularities; I do not know if one can construct a physical
model with these singularities.

If $\bar X_R$ can be smoothed, then the Hodge numbers
of the resulting manifolds will be different.
 \medpagebreak

\demo{Proof} We describe in details the case (2.9.1); the
others are similar. See also \S 6.
 $C_R$ is a $\bp $ bundle over the rational curve
$\phi _R(C_R)$, with fiber $f$ while $\mu _R (C_R)$ is a surface. The elliptic
 fibration $\pi _{S}: \pi ^{-1} (f)=S \to f$ is a rational elliptic surface 
with 
section (see [MV], [MS]), for each fiber $f$; the section is given by the
 intersection of $C_R$ with $S$.
 After the ``flop"  $S$ is a Del Pezzo surface $\mu _R(S)$, isomorphic to the 
blown 
up of $\Bbb P^2$ at $8$ points. 
 Each surface can be contracted to a point; actually all the surfaces
can be simultaneously contracted to a rational curve $\g _R$ (see 6.8), with a
 birational morphism $X_R \to \bar X _R$. Let $\rho: X \to \bar X_R$
denote the compositions of the two birational morphisms; from the explicit
construction of the flop it is clear that the elliptic fibration over $B_R$ is 
preserved and the following diagram is commutative 
 $$
\CD
X      @>\rho>>    \bar X_R \\
@V{\pi}VV           @VV{\bar \pi }V   \\
B        @>\phi_R>>   \bar B_R. \\ 
\endCD 
$$

Note that $ \codim \rho(D) \geq 2$, that is the image of the divisors 
contributing to the
superpotential is no longer a divisor. On the other hand $\bar X_R$ is singular
 along $\g _R$; these singularities are \it canonical \rm (like the 
singularities of the Weierstrass
model of $X$), [K], 1.5. $X_R$ is equisingular along $\g _R$: the 
singularity 
at each point of $\g _R$ of a transverse threefold
is exactly as in [MV,II].
 In fact, we can smooth $\bar X_R$ as in [MV,II]. 
\enddemo

The transitions among Fano threefolds with exceptional divisors
contributing to the superpotential appear in \S 7.

\medpagebreak

\definition{Example 5.3} \it (Where have all these divisors gone?) \rm
>From the Tables in \S 7, we can see  the sort of the other divisors
contributing to the superpotential after a birational contraction (5.1)-(5.2): 
some still contribute to the superpotential;
in some other cases the birational morphism $\phi_R$ becomes
a $\Bbb P^1$ (or conic bundle) fibration of $B_R$.

\smallpagebreak

 In example 4.4 ($B = S \times \bp$, with $S$ a rational elliptic surface) the 
birational
transformations $\phi_{R_\alpha}$ corresponding to the extremal ray $s_\alpha 
\times t$ 
contracts $B$ to $B_R=\bp \times S_1$, which is the unique Fano $3$-fold
with $b_2=10$ ($S_k$ is the Del Pezzo surface obtained by blowing up $\Bbb 
P^2$ at $9-k$ points; set $S=S_0$). 

We can perform $ 0 \leq k \leq 9$ contractions of non intersecting extremal 
rays and consider the induced elliptic fibration
$\pi _k : X \to  B_k= S_{9-k} \times \bp$. If $k>1$ there are $D_k$ smooth 
divisors
mapping to curves in $B_k$ and $\pi _k$ has no section.
  There is an infinite
number of divisors contributing to the superpotential,
whose image in $B_k$ is a singular divisor.
 I do not know at this moment if this can occur when there is
a section.

We can also contract the divisors $D_k$, as in 5.1 and 5.2 above.
A finite number of such divisors will still contribute to the superpotential, 
while infinitely many become singular
divisors with normal crossings.  $\diamondsuit$ \enddefinition

\vskip 0.5in

\S 6 Transitions of CY $4$-folds II:  \it Are these divisors ``exceptional"? 
\rm

\vskip 0.2in

We present some evidence that all the divisors
contributing to the superpotential (also in $M$-theory)
are ``exceptional", in the sense that are related to some birational 
transformation. They might not all the be exceptional,  in strict sense, as one 
can see in the
the example considered in [KV]. They show, that under certain hypothesis, if
$S \subset B$ is a rational surface and the ``general" fiber over a point in 
$S$ is a cycle of $N$ rational curves
with enhanced gauge group $SU(N)$, then each of the $N$
irreducible component of $SU(N)$  contribute to the superpotential. However 
only $N-1$ of them are
``exceptional" divisors. In this case the birational morphism
is the contraction to the Weierstrass model (1.2).
It should be pointed out that there exists a relation
among these $N$ divisors ($N-1$ are ``independent") [KV].

\medpagebreak

 If the normal bundle is negative (Grauert)  a
contraction is possible, at least in the analytic category.
We would like this contraction to be  projective
and to  describe the singularities which might occur.
 In the case of $F$-theory we would like also
to preserve the elliptic structure.

\medpagebreak

Our approach is to consider the pair $(X,D)$, where
$D$ is a divisor contributing to the superpotential
and exploit once more the fact that this divisor cannot be
\it nef \rm (3.1). 
 We will need some more general definition that the ones in \S 2.

 The reader should probably start from the second part of
this section (``the general case") and use the first one as
a reference.

\vskip 0.2in

{\bf 6.a  Log minimal models}

\bigpagebreak

There are several version of the log minimal model
program; we follow [KMM], as it seems at the moment to
be the best suited for our applications.

$\pxb$ is any proper morphism between varieties; later
we will apply
the general machinery to the case of the elliptic Calabi-Yau.

\definition{Definition 6.1}$\nxb \subset \Bbb R ^m$ is the closed convex cone
 generated (over $\Bbb R _{\geq 0}$) by the effective cycles of (complex)
dimension 1, mod. numerical equivalence.
\enddefinition 

\definition{Definition 6.2} $D$ is $\pi$-\it nef \rm if $D \cdot \g \geq 0$,
for all the curves $\g \in \nxb$.
\enddefinition 

A \it relative \rm version of Kleiman's criterion says that
the cone of $\pi$-nef divisors (which is the closure of the $\pi$-ample cone)
and $\nxb$ are dual cones.
 The duality is again given by the intersection pairing (2.0).

In what follows we will have to consider  singular varieties;
a crucial point in the (log)-minimal model program is
the existence of a ``reasonable" intersection pairing between
complex curves (with values in $\Bbb Q$)
and complex subvarieties of codimension $1$
(\it Weyl divisors \rm). This motivate the following:
 
\definition{Definition 6.3} A variety has $\Bbb Q$-factorial singularities if 
for any $D$  Weil divisor, there 
exists and integer $r$ such that $rD$ is a line bundle.
($D$ is also called $\Bbb Q$-Cartier divisor.)
\enddefinition

Unless noted otherwise all the varieties are assumed to be normal and $\Bbb 
Q$-factorial. We will also consider
Weyl divisors with rational coefficients.

  Below $\Cal D $ is such a divisor: $ \Cal D=\sum a_i L_i$, with $ L_i$ 
distinct 
complex subvarieties of codimension $1$ (Weyl divisors) and $a_i \in \Bbb Q$, 
$0 \leq a_i  < 1$. 
$\cup L_i$ is called \it support \rm of $\Cal D$.

We write  $\Cal D\equiv \Cal D'$ if some multiple
of $\Cal D$ and $\Cal D'$ are equivalent as line bundles.

\definition{Definition 6.4} The pair $(X, \Cal D)$ (as above)
has at worse 
\it log-terminal \rm (\it log canonical\rm) singularities
if there exists a resolution of the singularities $f: Y \to X$
such that the union of the exceptional divisor
and the inverse image of $ \cup L_i$ is a divisor with normal
crossings and
$$K_Y \equiv  f ^*(K_X + \Cal D) + \sum b_k M_k ,\text{ such that } \ b_k > -1 
 \ (\text{resp.} \geq  -1)  , \ \forall k.$$
\enddefinition
(The definition does depend on the choice of $f$ and $Y$.)
If $\Cal D =0$, and $b _k >0 $ ($b_k \geq 0$)
then the singularities are called
\it terminal \rm and \it canonical \rm respectively.

If $\dim(X)=2$ and $X$ the singularities are at worse terminal,
then $X$ is smooth, the canonical singularities are the rational double points.

The following is a generalized version of the contraction
theorem (2.6):

\proclaim{Theorem 6.5 (Contraction morphism)} Let $\pxb$ be a morphisms 
between varieties.
If $(X, \Cal D)$ has log-terminal singularities and  $K_X + \Cal D$ is not 
$\pi$-nef (that is $(K_X + \Cal D) \cdot R <0$, for
some extremal ray $R \in \nxb$), then there exists a morphism 
that  $\psi_R : X \to 
Z$, contracting all the curves in the numerically equivalence (homology) class 
of $[R]$ such that the following diagram is commutative.
$$
\alignat3
& {X} &&  @>\psi _R>> & \ \ \ {Z}\\
& _{\pi } \searrow && & \swarrow _{ \pi _+}\\
& &&\ \ B &&.
\endalignat
$$

$Z$ is a normal variety and $\dim \nxb > \dim \overline {NE(Z/B)}$.
\endproclaim

\demo{Proof} For a proof and various reference, see
for example [KMM], Theorems 3.1.1,  3.2.1, 4.1.1 and 4.2.1.\qed
\enddemo

(*) We assume also that some line bundle multiple of
$K_X + \Cal D$ has a section (i.e. the Kodaira dimension
of $K_X + \Cal D$ is non negative).
 This is the case in our applications, where $K_X \sim \Cal O_X$ and a
 multiple of $\Cal D$ is an effective divisor contributing to the
 superpotential.

In this case, the contraction morphism in (6.5) is
birational.
\medpagebreak

The \it Log Minimal Model Conjecture \rm says that  there exists
a birational map $\rho: X \dasharrow \bar X$ and
a morphism $\bar \pi : \bar X \to B$
such that  $K_{\bar X} + \bar {\Cal D}$ is $\bar \pi$-nef
and the following diagram is commutative:
$$
\alignat3
& {X} &&  @>\rho>> & \ \ \ {\bar X}\\
& _{\pi } \searrow && & \swarrow _{ \bar \pi }\\
& &&\ \ B &&.
\endalignat
$$
 Here $\bar {\Cal D} = \rho (\Cal D)$ and $(\bar X, \bar {\d})$ is the \it 
log-minimal model \rm.
 
\smallpagebreak

 The problem is that when the  contraction in (6.5)
is not divisorial (that is 
the exceptional locus is not a divisor), it is not possible to define an 
intersection product which is compatible with our structure (6.3). If so, we 
would in fact have a contradiction:

 $0 > (K_X + \Cal D) \cdot R = {\psi _R} ^*
(K_Y +  \bar {\d}) \cdot R = (K_Z +  \bar {\d}) \cdot \psi _R (\d)= 0.$

 In this case there is the following:

\proclaim{Conjecture 6.6} There is another birational transformation ( 
``log-flip"), which   is an isomorphism
outside a set of codimension greater than $2$ (an isomorphism
in codimension $1$):
\medpagebreak

\noindent\phantom{it's alongwaytoTippe}
\epsfxsize=2in\epsfbox{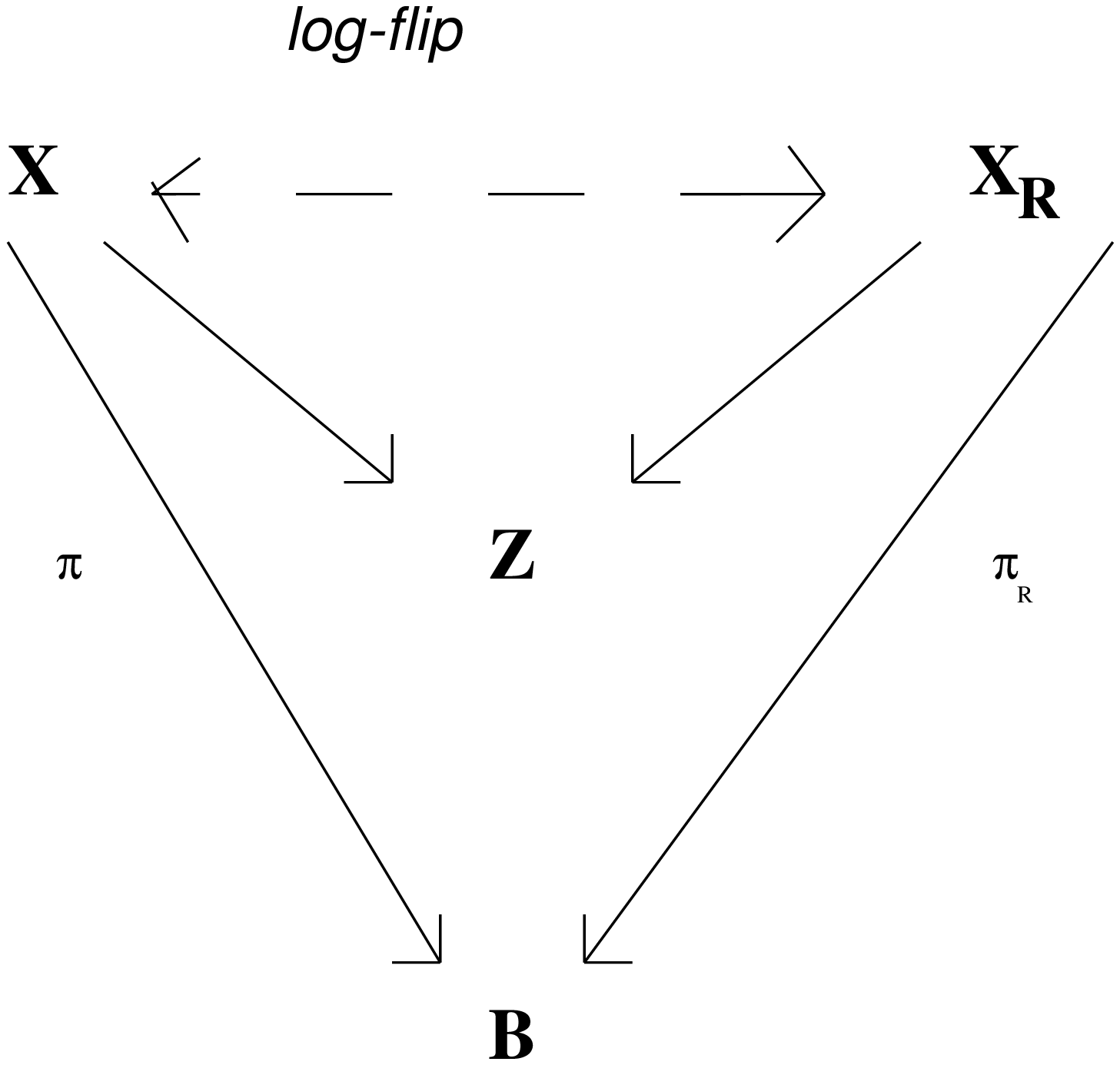}

\medpagebreak

such that $X_R$ has log-terminal ($\Bbb Q$-factorial) singularities and 

\noindent $\mu (K_{X_R} + \mu (\d)) \cdot {R_+} >0$,
for all the curves $A$ contracted by $\psi _+$.

The number  consecutive 
such log-flips is always a finite.
\endproclaim

The Log Minimal Model Conjecture is a theorem if
$\dim(X) \leq 3$ (see for example [K.etal]) and
it has been be worked out in various special examples, among which the
ones considered in \S 5 (which we review in (6.8) below) and when the 
techniques of toric geometry can be
applied [O], [Mt].

\vskip 0.2in

{6.b \bf Transition II: the general case}

\bigpagebreak

Now let $X$ be a smooth Calabi-Yau $4$-fold and  $D$ a divisor contributing to 
the superpotential; then $D$ is not nef, 
that is, there is an effective (complex) curve $R$
such that $D \cdot R = (K_X + D) \cdot R<0$ (3.1). The idea is to use the 
contraction
morphism in (6.5) 

\smallpagebreak

We consider the pair $(X, \d)$, where $\d = rD$,
for some $0<r <1, \ \in \Bbb Q$:
$X$ is smooth, so we can take as $f$ in (6.4)  the identity map and verify 
that 
the pair has log terminal singularities (this is true also if $D$ has normal 
crossings singularities).

\medpagebreak

 If the log-minimal model conjecture holds, then
the following  conjecture is true:

\proclaim{Conjecture 6.7 } Let $X$ be a Calabi-Yau
fourfold and $D$ a divisor contributing to the superpotential.
Then  there exists a birational transformation $\rho: X \to \bar X$,
 with canonical singularities (the same singularities of
the Weierstrass model) and $\rho (D)$ is a nef effective divisor.
\endproclaim

\demo{Proof} Start with $(X, \Cal D)$ as above.
If the log-minimal model conjecture holds, then
$\rho$ is a composition of contraction morphisms (6.5)
and log-flips (6.6).
If $\nu : X \to X'$ is either the contraction morphism in (6.5), or the
``log-flip" in (6.6)
 then $X'$ has canonical singularities ([Ka],  1.5 ); these are the same
singularities of our Weierstrass models (1.1)).
 Then $K_X' \sim \Cal O_{X'}$ and $K_{X'} + \nu (\Cal D) \sim \nu (\Cal D)$. 

Note that $\nu (\Cal D)$ is well defined and that these
log-flips are ``flops" (because the canonical divisor
is trivial). \qed
\enddemo

 (5.1) and (5.2) are particular cases of this general set up.

\definition{Example 6.8} Let $\pxb$ an elliptic fibration between smooth 
varieties.
Assume that $X$ is equal to the  smooth ``general" Weierstrass
model over $B$ and that $D$ is a divisor contributing to
the superpotential. Then $D=\pi ^*(C)$, for some
smooth divisor on $B$.

Now let us consider the induced elliptic fibration 
$\epsilon: X \to B_R$ and $\overline {\NE(X/B_R})$.
 This two dimensional cone is generated by a fiber
$\g$ of the fibration $\pi$ and the extremal ray $R$ in $X$
(more precisely, the isomorphic image in the section $T \subset
X$ of the extremal ray $R$ in $B$):
$D\cdot \g =0$, while $D \cdot \g <0$. 

Then there exists
a contraction morphism 
$\psi_R : X \to Z$ (6.5) contracting the curves in the homology
of class $[R]$; this contraction cannot be divisorial
(it comes from a contraction from the lower dimensional $B$). In each of the 
cases considered
the flop $\mu _R:  X \dasharrow X_R$
exists [Mt]. Let $\pi _R$ be the induced elliptic fibration.
We now concentrate on the case (2.9.1); the others are similar.

 After the ``flop" the relative cone $\overline \NE(X_R/B_R)$
is still two dimensional and it is generated by
the image of the fiber of $\pi _R$, which we still denote by $\g$,
and $R_{+}$, a fiber of $ X _R \to Z$.
 It is easy to verify that $\mu(D) \cdot R_{+} >0$, while
$\mu(D) \cdot \g <0$.
In this case the contraction morphism corresponding to $\g$
is divisorial
and the divisor $\Cal D$ is the exceptional divisor.
\medpagebreak

\noindent \phantom{itsa longwayto} \epsfxsize=3in\epsfbox{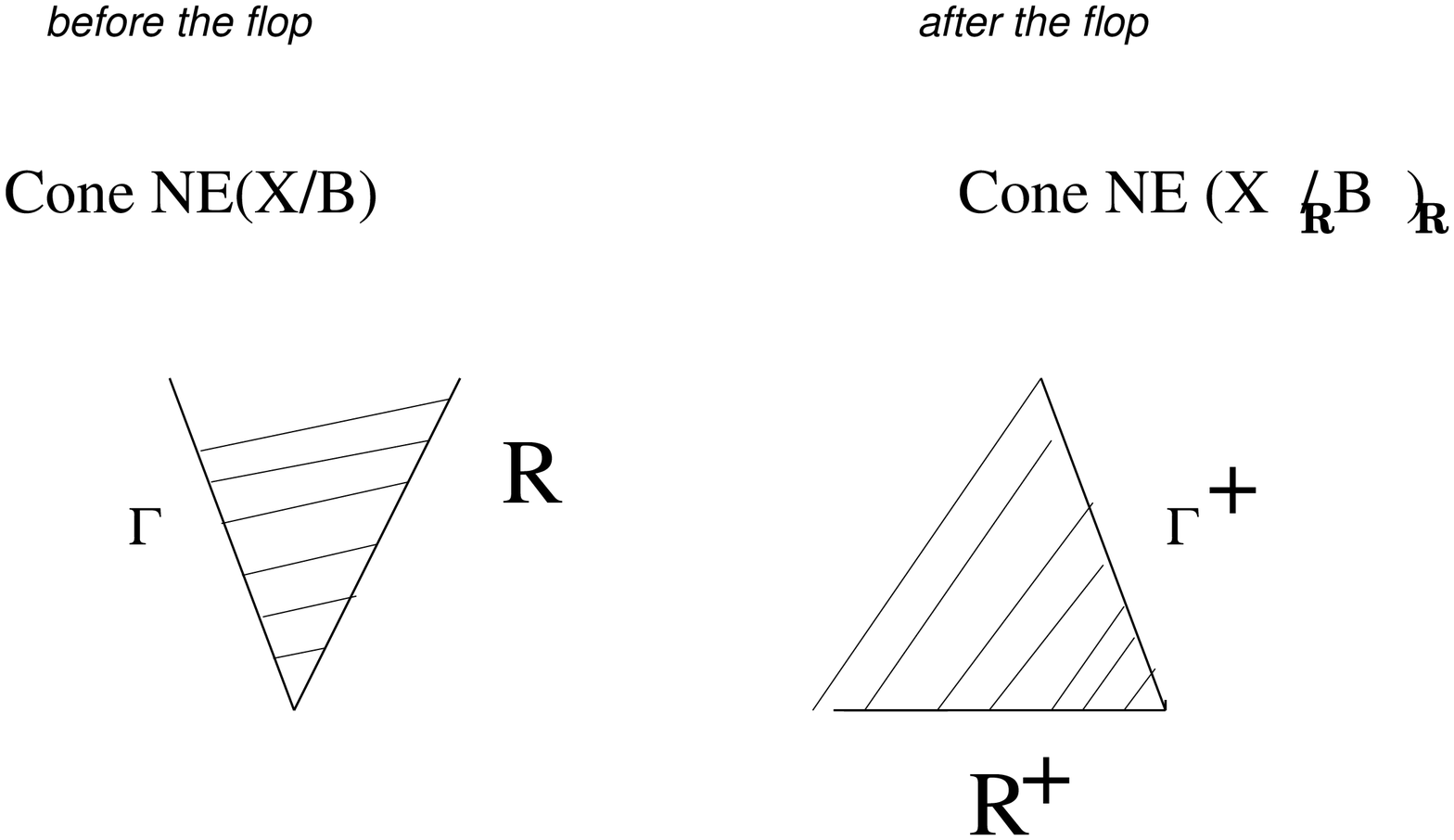}

\medpagebreak

If $\pxb$ has more than one section (the rank of the Mordell-Weyl group is 
positive)
we first have to
perform ``flops" along the sections (as in [MV,II]). $\diamondsuit$
 \enddefinition

\smallpagebreak

\newpage

\S 7 Tables for ``general" elliptic CY with basis Fano $3$-folds.

\bigpagebreak

In the following, $\pxb$ is an elliptic Calabi-Yau $4$-fold and $B$ is
a Fano $3$-fold. In Tables 7.1-7.6
we follow
 Iskovskih-Mori-Mukai's list of Fano $3$-folds: the
threefolds are subdivided by their second betti number, $1 \leq b_2= h^{1,1} 
\leq 10$, which is
also the dimension of the Mori (and K\"{a}hler) cone of $B$.

 We gather various information
about the Fano threefolds and the
``general" elliptic Calabi-Yau $4$-folds fibred over them.

We use the criteria developed in \S 1 and \S 3 to determine the divisors of 
the form $D= \pi ^*(C)$
which contribute to the superpotential on $X$ (4.5.4).
If $X=W$, the smooth Weierstrass model (1.0),
these divisors are all the divisors contributing to the superpotential
 (1.2). The divisors $C$  determine
a birational contraction $B \to B_R$. We identify
$B_R$ when it is another Fano (4.5.1).
 If $X=W$  we also compute the topological euler characteristic of $X$.
In particular,

\medpagebreak

$\bullet$ The first number in the table corresponds to the one assigned in [MM]
to each $3$-fold with a given $b_2= h^{1,1}$.
If $\pi: W=X \to B$ is the smooth general Weierstrass model (1.0) (there is
only one section of $\pi$), then $h^{1,1}(X)= h^{1,1}(B) +1$. 

\medpagebreak

$\bullet$  The second column says whether $B$ is toric: a list of toric 3 fold
and the related superpotential appears in [KLRY] and [Moh], if the $3$-fold
is toric $\Cal F_k$ is the symbol used in [B], [O] and [Moh]; many examples
are also in [My].

\medpagebreak

$\bullet$  The third column is about the divisors contributing to the 
superpotential
as in 4.5.1.  If $D_j= \pi ^*(C)$, then $C$ is a divisor of type
(2.9.j), $1 \leq j \leq 4$; if there are 2 different
divisors of the same type (2.9.j) we will denote them as $D_j^1, \ D_j^2$.

\smallpagebreak

 If the same divisors is exceptional for 2 different contractions (as in 
4.5.1)  we simply write it twice 
(this is the case of no. 3,  $h^{1,1} =3$).

 If $\phi : B \to B_R$ is the contraction of $D_j$
then $h ^{1,1} (B_R) = h^{1,1}(B)-1$. If $X=W$ these birational 
transformations  are
``promoted" to birational transitions of the Calabi-Yau
fourfold $X$ (see \S 5 and \S 6.b).

\smallpagebreak

  If $B_R$ is  Fano,
$D_j(\ell)$ means that $B_R$ is the Fano $3$-fold with number
$\ell$ in the Mori-Mukai classification [MM] of threefolds with $h^{1,1}= 
h^{1,1}(B)-1$.

\medpagebreak

$\bullet$  The fourth column lists the $\bp$-fibrations
 (denoted by $p_i: B \to S$)
and the conic bundles (denoted by $c_i: B \to S$): this is relevant
from the point of view of heterotic theory.

\medpagebreak

$\bullet$ The fifth column is $12c_1(B) \cdot c_2(B) + 360 c_1 ^3(B) $,
 which is the Euler characteristic of the smooth
Weierstrass model (if any) over $B$ [SVW].
 By the Riemann-Roch theorem for threefolds [H],
 $12c_1(B) \cdot c_2(B)=288 \chi (\Cal O_B) = 288$ ($B$ is uniruled).

\medpagebreak

$\bullet$  We use a rather    crude  (but readily available [I], [MM2]) 
criterion
 to determine whether there exists
a smooth Weierstrass Calabi-Yau model  over $B$ (1.0), namely we
require  $-K_B$ to be very ample. On the other hand, most Fano satisfy this 
criterion: we write ``no" in the last column
if $-K_B$ is not very ample. Otherwise $X=W$, its smooth Weierstrass model; in 
this case, we see from the list that: $$\chi(X)=144(17+5 \ell), \ 0 \leq \ell 
\leq 25, \ \ell=28, 29.$$

\medpagebreak

$\bullet$ {\bf Table 7.7 } is the flow chart of transition among the
Fano threefolds corresponding to divisors contributing
to the superpotential (as in 4.5.1). These 
are also ``promoted" to transition among
Calabi-Yau threefolds (as in \S 5 and \S 6.b).

 The columns correspond the the values of
$h^{1,1}(B)$, starting from 5 on the left
and ending with 1 on the right.
\medpagebreak

 The thick lines represent a contraction of a divisor
($\sim \bp \times \bp$) to a point (2.9.2),
while the others represent a contraction of a rational
ruled surface (2.9.1).

\vskip 0.2in

\centerline{ \bf THE TABLES} 

\bigpagebreak

{\bf Table 7.1: $h^{1,1}(B)=1$}

\medpagebreak

 Iskovskih [I] classified all such varieties: the following occur
in the flow chart (Table 7.7), together with $\Bbb P^3$:

\noindent $\ast$ $Q \subset \Bbb P^4$, a smooth quadric surface;

\smallpagebreak

\noindent $\ast$ $V_3 \subset \Bbb P^3$, a smooth cubic surface;

 \smallpagebreak

\noindent $\ast$ $V_4 \subset \Bbb P^5$, a complete intersection of two
quadrics;

\smallpagebreak

\noindent $\ast$ $V_5 \subset \Bbb P^9$ is a complete intersection of a linear 
subspace $\Bbb P^6 \subset \Bbb P^9$ and the Grassmann variety of $\Gr (1,4)$
embedded in $\Bbb P^9$ by the Plucker embedding.

\medpagebreak

$\bullet$ The only Fano toric variety with $h^{1,1}(B)=1$ is $\Bbb P^3$.

\medpagebreak

$\bullet$ No divisor contribute to the superpotential (4.1) and there is no
$\bp$-fibration.

\medpagebreak

 $\bullet$ All  these $3$-folds have $-K_B$ very ample
with the following exceptions  (see [I], vol 12, table 6.5
or also [Mur]):

\smallpagebreak

\noindent $\ast$ the double cover of $\Bbb P^3$ with branch locus a sextic,

\noindent $\ast$ the double cover of a quadric in $\Bbb P^4$ branched over the 
intersection of the quadric and a quartic

\noindent $\ast$ $V_1$ (i.e. the double cover of the cone over the Veronese)

\noindent $\ast$ $V_2$ (i.e. the double cover of $\Bbb P^3$ with
quartic ramification).

\vskip 0.5in

The next page contains {\bf Table 7.2} for $h^{1,1}(B)=2$.

\newpage

\vbox{\tabskip=0pt \offinterlineskip
\halign to 321pt
{\vrule \hfil\strut#\hfil
&\vrule  \hfil#\hfil&\vrule\hfil#\hfil&\vrule\hfil#\hfil
&\vrule\hfil#\hfil\vrule&\vrule\hfil#\hfil\vrule&\vrule\hfil#\hfil\vrule
\cr
\tablerule \ no.    & \ tor.   & contr. to the superpo.  & 
  fibrations &  \ $\chi (X)$& \ v.a.  \cr
\tablerule \ 1    & \  no   & \  none   & 
 \quad none  \quad \quad & \   & \ no \cr
\tablerule \  2    & \ no    & \  none   & 
\quad $c:B \to \Bbb P^2$ &  \ & no \cr
\tablerule \ 3    & \  no   & \  none   & 
 \quad none  \quad \quad & \    & no \cr
\tablerule \ 4    & \  no   & \  none   & 
 \quad none  \quad \quad &\ 3888  & \cr
\tablerule \ 5    & \  no   & \  none   & 
 \quad none  \quad \quad &\ 4608    & \cr
\tablerule \ 6    & \  no   & \  none   & 
$c_1:B \to \Bbb P^2$, $c_2:B \to \Bbb P^2$ &\ 4608  & \cr
\tablerule \ 7    & \  no   & \  none   & 
\quad none  \quad \quad  &\ 5328 &   \cr
\tablerule \ 8.a    & \  no  &  $D_4$   & 
\quad $c:B \to \Bbb P^2$  \quad \quad  &\ 5328 & \cr
\tablerule \ 8.b    & \  no  & \ none   & 
\quad $c:B \to \Bbb P^2$  \quad \quad  &\ 5328 &   \cr
\tablerule \ 9    & \  no  & \ none   & 
\quad $c:B \to \Bbb P^2$  \quad \quad  &\ 6048  & \cr
\tablerule \ 10    & \  no   & \  none   & 
 \quad none  \quad \quad &\ 6048 & \cr
\tablerule \ 11    & \  no   & \  $D_1$ ($V_3$)  & 
$c:B \to \Bbb P^2$ & \  6768 & \cr
\tablerule \ 12    & \  no   & \  none   & 
 \quad none  & \ 7488 & \cr
\tablerule \ 13    & \  no   & \  none   & 
 $c:B \to \Bbb P^2$  & \ 7488   & \cr
\tablerule \ 14    & \  no   & \  none   & 
 \quad none  \quad \quad &\ 7488  &  \cr
\tablerule \ 15.a    & \  no   & \  $D_4$   & 
 \quad none  \quad \quad & \ 8208   & \cr
\tablerule \ 15.b    & \  no   & \  none   & 
 \quad none  \quad \quad &\ 8208 & \cr
\tablerule \ 16    & \  no   & \  $D_1$ ($V_4$) & 
$c:B \to \Bbb P^2$ & \  8208 & \cr
\tablerule \ 17    & \  no   & \  none   & 
 \quad none  \quad \quad & \ 8928   &  \cr
\tablerule \ 18    & \  no  & \ none   & 
\quad $c:B \to \Bbb P^2$  \quad \quad  &\ 8928 &  \cr
\tablerule \ 19    & \  no   & \  $D_1$ ($V_4$)   & 
 \quad none  \quad \quad & \ 9648  &  \cr
\tablerule \ 20    & \  no  & \ $D_1$ ($V_5$)  & 
\quad $c:B \to \Bbb P^2$  \quad \quad  &\ 9648 &  \cr
\tablerule \ 21    & \  no   & \  $D^1_1$ ($Q$), \quad $D^2 _1$   & 
 \quad none  \quad \quad & \ 10368  &  \cr
\tablerule \ 22    & \  no   & \  $D^1_1$ ($\Bbb P^3$), \quad $D^2 _1$ ($V_5$) 
 & 
 \quad none  \quad \quad & \ 11088  &  \cr
\tablerule \ 23.a    & \  no   & \  $D_4$   & 
 \quad none  \quad \quad & \ 11088  &  \cr
\tablerule \ 23.b   & \  no   & \  none   & 
 \quad none  \quad \quad & \ 11088  &  \cr
\tablerule \ 24    & \  no   & \  none   & 
$c:B \to \Bbb P^2$, $p:B \to \Bbb P^2$ &\ 11088  & \cr
\tablerule \ 25    & \  no   & \  none   & 
 \quad none  \quad \quad &\ 11808  &  \cr
\tablerule \ 26    & \  no   & \  $D^1_1$ ($Q$), \quad $D^2 _1$ ($V_5$)   & 
 \quad none  \quad \quad & \ 12528  &  \cr
\tablerule \ 27    & \  no   & \  $D_1$ ($\Bbb P^3$)  & 
$p: B \to \Bbb P^2$ & 13968 & \cr
\tablerule \ 28    & \  no   & \   \quad $D _3$   & 
 \quad none  \quad \quad & \ 14688   & \cr
\tablerule \ 29    & \  no   & \   \quad $D _1$ ($Q$)  & 
 \quad none  \quad \quad & \ 14688  &  \cr
\tablerule \ 30    & \  no   & \  $D_1$ ($\Bbb P^3$), \quad $D_2$   & 
 \quad none  \quad \quad & \ 16848  &  \cr
\tablerule \ 31    & \  no  & \ $D_1$ ($Q$) & 
\quad $p: B \to \Bbb P^2$  \quad \quad  &\ 16848 &  \cr
\tablerule \ 32    & \  no   & \  none   & 
$p_1: B \to \Bbb P^2$, $p_2: B \to \Bbb P^2$ &\ 17568 &   \cr
\tablerule \ 33    & \ $\F _5$   & \  $D_1$ ($\Bbb P^3$)  & 
 \quad none  \quad \quad & \ 19728  &  \cr
\tablerule \  34    & $\F _2$    & \  none   & 
\quad $p: B \to \Bbb P^2$ &  \ 19728 & \cr
\tablerule \ 35    & $\F _3$  & \ $D_2$ ($\Bbb P^3$) & 
\quad $p: B \to \Bbb P^2$  \quad \quad  &\ 20448 & \cr
\tablerule \  36    &  $\F _4$   & \  none   & 
\quad $p: B \to \Bbb P^2$ &  \ 22608 & \cr
\tablerule }}
\newpage

{\bf Table 7.3: $h^{1,1}(B)=3$}

\vskip 0.3in

{\vbox{\tabskip=0pt \offinterlineskip
\halign to 365pt
{\vrule \hfil\strut#\hfil
&\vrule  \hfil#\hfil&\vrule\hfil#\hfil&\vrule\hfil#\hfil
&\vrule\hfil#\hfil&\vrule\hfil#\hfil&\vrule\hfil#\hfil&\vrule\hfil#\hfil\vrule
\cr
\tablerule \ no.    & \ tor.   & contr. to the superp.& 
  fibrations & $\chi (X)$ & \cr
\tablerule \ 1    & \  no   & \  none   & 
 $ c_1:  \bp \times \bp ,   c_2:   \bp \times \bp ,  
c_3:  \bp \times \bp$  & \ 4608  & \cr
\tablerule \  2    & \ no    & $D_1, D_1$   & 
\ $c:B \to \bp \times \bp$ &  \ 5328 &  \cr
\tablerule \ 3    & \  no   & \  none   & 
 \  $c:B \to \bp \times \bp$ & \ 6768    &  \cr
\tablerule \ 4    & \  no   & \ $D_1$ (18)  & 
 \ $c_1 : B \to \Bbb P^2,  \ c_2 : B \to  \bp \times \bp$ &\ 6768   &  \cr
\tablerule \ 5    & \  no   & $D^1 _1 (34), \ D^2 _1, \ D^2 _1$   & 
 \ none  \ &\ 7488    &  \cr
\tablerule \ 6    & \  no   & \  $D_1$ (33)   & 
\ $c:B \to \bp \times \bp$ &\ 8202   &  \cr
\tablerule \ 7    & \  no   & \  none   & 
\ none    &\ 8928    &  \cr
\tablerule \ 8    & \  no  & \ $D^1 _1 (24), \ D^2 _1$ (34) & 
\quad $c:B \to \bp \times \bp$  &\ 8928    &  \cr
\tablerule \ 9    & \  no  & \ $D^1_3, \ D ^2 _3$  & 
\ none &\ 9648   &  \cr
\tablerule \ 10    & \  no   &  $D^1 _1 (29), \ D^2 _1 $  & 
 \  none &\ 9648  &  \cr
\tablerule \ 11    & \  no   & \  $D^1_1 (25), \ D^2 _1$ (34)  & 
 \ none & \  10368 &  \cr
\tablerule \ 12    & \  no   & \ $D^1_1 (27), \ D^2_1 (33), \ D^3_1 (34)$   & 
 \ none  & \ 110368 &  \cr
\tablerule \ 13    & \  no   & \ $D^1_1 (32), \ D^2 _1, \ D^3 _1$   & 
\ none  & \ 11088    &  \cr
\tablerule \ 14    & \  no   & \  $D_2 (28), \ D_3$   & 
 \ none &\ 11808    &  \cr
\tablerule \ 15    & \  no   & \ $D^1_1 (29), \ D^2_1 (31), \ D^3_1 (34)$   & 
 \ none   &\ 11808  &  \cr
\tablerule \ 16    & \  no   &\ $D^1_1 (27), \ D^2_1 (32), \ D^3_1 (35)$ & 
\ none  & \ 12528 &  \cr
\tablerule \ 17    & \  no   & \ \ $D^1_1 (34), \ D^2_1$   & 
 \ $ p: B \to \bpt$ & \ 13248    &  \cr
\tablerule \ 18    & \  no  & \ $D^1_1 (29), \ D^2_1 (30), \ D^3 _1 (33)$    & 
\ none  &\ 13248   &  \cr
\tablerule \ 19    & \  no   &   $D^1_1 (35), \  D^2_1, \  D^3 _1, \ D^4_1$  & 
 \ none  & \ 13968    &  \cr
\tablerule \ 20    & \  no  & \ $D^1_1 (31), \  D^2_1 (32), \  D^3 _1$  & 
\ none &\ 13968   &  \cr
\tablerule \ 21    & \  no   & \  $D^1_1 (34), \  D^2_1, \  D^3 _1$  & 
 \ none   & \ 13968    &  \cr
\tablerule \ 22    & \  no   & \  $D^1_1 (34), \  D^2_1 (36), \  D _3$   & 
 \ none   & \ 14688    &  \cr
\tablerule \ 23   & \  no   & \ $D^1_1 (30), \  D^2_1 (31), \  D^3 _1 (35) $  
& 
 \ none  & \ 15408    &  \cr
\tablerule \ 24    & \  no   & $D^1_1 (32), \  D^2_1 (34)$  & 
$p:B \to \Bbb F _1$ &\ 15408  &  \cr
\tablerule \ 25    & $\F _8$   & $D^1_1 (33), \  D^2_1$ & 
 \ $p: B \to \bpt$ &\ 16128   &  \cr
\tablerule \ 26    & $\F_ {12}$   & \  $D^1_1 (33)\ D^2 _1 (34), D_2 (35)$   & 
 \ none  & \ 16848    &  \cr
\tablerule \ 27    &  $\F_6$  & \  none  & 
$p_1:  \bpt, p_2:  \bpt, p_3:  \bpt$ & 17568 &  \cr
\tablerule \ 28    & $\F _9$   & \    $D _1$ (34)  & 
 \ $p_1: B \to \bpt, p_2: B \to \bpt$  & \ 17568    &  \cr
\tablerule \ 29    & $\F_{ 11}$  & \   $D^1_1 (35), \ D^2 _1 (36), \ D_3$  & 
 \ none  & \ 18288    &  \cr
\tablerule \ 30    & $\F _{10}$   & \  $D^1_1 (33), \ D^2_1 (35)$   & 
 \ $p: B \to \Bbb F_1$ & \  18288   &  \cr
\tablerule \ 31    & $\F _7$  & \ $D_1, \ D_1$  & 
\  $p: B \to \bpt$   &\ 19008 & \cr
\tablerule }}

\newpage

$S_k$ denotes $\Bbb P^2$ blown up at $9-k$ points; for example
$\F _{13} = \bp \times S_7$.

\bigpagebreak

{\bf Table 7.4:} $h^{1,1}(B)=4$.

\bigpagebreak

{\vbox{\tabskip=0pt \offinterlineskip
\halign to 312pt
{\vrule \hfil\strut#\hfil
&\vrule  \hfil#\hfil&\vrule\hfil#\hfil&\vrule\hfil#\hfil
&\vrule\hfil#\hfil\vrule&\vrule\hfil#\hfil\vrule&\vrule\hfil#\hfil\vrule
\cr
\tablerule \ no.    & \ tor.   & contr. to the superp.& 
  fibrations & $\chi(X)$ \cr
\tablerule \ 1    & \  no   & \  none   & 
 \ none  & \ 8928    \cr
\tablerule \  2    & \ no    & $D^1_1 , D^1_1, D^2_1, D^2_1 $   & 
\ none &  \ 10368 \cr
\tablerule \ 3    & \  no   &  $D^1_1 (17), D^2_1 (27), D^3_1 (28), D^4_1 (28) 
$  & 
 \  none & \ 11088    \cr
\tablerule \ 4    & \  no   & \ $D^1_1 (18), D^2_1 (18), D^3_1 (19), D^4_1 
(30), D^5_1 (30) $    & 
 \ none &\ 11808   \cr
\tablerule \ 5    & \  no    & $D^1_1 (21), D^2_1 (28), D^3_1 (31), D^4_1, 
D^5_1 $    & 
 \ none &\ 11808   \cr
\tablerule \ 6    & \  no   &  $D^1_1 (25), D^2_1 (25), D^3_1 (25), D^4_1 (27) 
$  & 
 \  none & \ 12528    \cr
\tablerule \ 7    & \  no   &  $D^1_1 (24), D^2_1 (24), D^3_1 (28), D^4_1 (28) 
$  & 
 \  none & \ 13248    \cr
\tablerule \ 8    & \  no  &  $D^1_1 (27), D^2_1 (31), D^3_1 (31), D^4_1 (31) 
$  & 
 \  none & \ 13968    \cr
\tablerule \ 9    & $\F _{15}$   &  $D^1_1 (25), D^2_1 (26), D^3_1 (28), D^4_1 
(30) $  & 
 \  none & \ 14688    \cr
\tablerule \ 10    & $\F _{13}$  &  $D^1_1 (27), D^2_1 (28), D^3_1 (28)$  & 
 \ $p: B \to S_7$ & \ 15408  \cr
\tablerule \ 11    & \   $\F _{14}$   & $D^1_1 (28), D^2_1 (31), D^3_1, D^4_1 
$  & 
 \  none & \ 16128     \cr
\tablerule \ 12    & \  $\F _{16}$   &\ $D^1_1 (30), \ D^2_1 (30), \ D^3_1, 
D^4_1$   & 
 \ none  & \ 16848 \cr
\tablerule }}

\vskip 0.3in

{\bf Table 7.5:} $h^{1,1}(B)=5$.

\medpagebreak

{\vbox{\tabskip=0pt \offinterlineskip
\halign to 367pt
{\vrule \hfil\strut#\hfil
&\vrule  \hfil#\hfil&\vrule\hfil#\hfil&\vrule\hfil#\hfil
&\vrule\hfil#\hfil\vrule&\vrule\hfil#\hfil\vrule&\vrule\hfil#\hfil\vrule
\cr
\tablerule \ no.    & \ tor.   & contr. to the superp.& 
  fibrations & $\chi(X)$ \cr
\tablerule \ 1    & \  no   & $D^i_1 (4),  i=1,2, 3,  D^i_1 (12), i=4,5 , 6,  
D^7_1 $  & 
 \ none  & \ 10368    \cr
\tablerule \  2    & $\F _{18}$     & $D^i_1 (9),  i=1, 2, D^i_1 (11), 
\  i=3, 4, \ D^5_1 (12), D^i_1, \ i=6,7$  & 
\ none &  \ 13248 \cr
\tablerule \ 3    & $\F _{17}$    &  $D^i_1 (10), \ i=1, \cdots , 6 $  & 
 \ $p: B \to S_6$ & \ 13248  \cr
\tablerule }}

\vskip 0.3in

{\bf Table 7.6:} $h^{1,1}(B) \geq 6$.

\medpagebreak

$B= \bp \times S_{k}$, with $1 \leq k \leq 5$; $h^{1,1} (B)= 11 -k$.
None of these threefolds is toric; the extremal contractions
are induced by the blow ups: $S_k \to S_{k+1}$.

\bigpagebreak

{\vbox{\tabskip=0pt \offinterlineskip
\halign to 243pt
{\vrule \hfil\strut#\hfil
&\vrule  \hfil#\hfil&\vrule\hfil#\hfil&\vrule\hfil#\hfil
&\vrule\hfil#\hfil\vrule&\vrule\hfil#\hfil\vrule&\vrule\hfil#\hfil\vrule
\cr
\tablerule \ $h^{1,1}(B)$    &  \ contr. to the superp  & fibrations & 
  $\chi(X)$  &  \ v.a. \cr
\tablerule \ 6    & $D^i_1, \ i=1, \cdots , 10 $   & \  $ p: B \to S_5$   & 
 13248  &     \cr
\tablerule \ 7    & $D^i_1, \ i=1, \cdots , 16 $   & \ $ p: B \to S_4$   & 
 15408  &     \cr
\tablerule \ 8    & $D^i_1, \ i=1, \cdots , 27 $   & \ $ p: B \to S_3$   & 
 17568   &      \cr
\tablerule \ 9    & $D^i_1, \ i=1, \cdots , 56 $  & \  $ p: B \to S_2$   & 
   &  no   \cr
\tablerule \ 10    & $D^i_1, \ i=1, \cdots , 240 $   & \ $ p: B \to S_1$   & 
   &  no    \cr
\tablerule }}

\newpage

\noindent \epsfxsize=5in\epsfbox{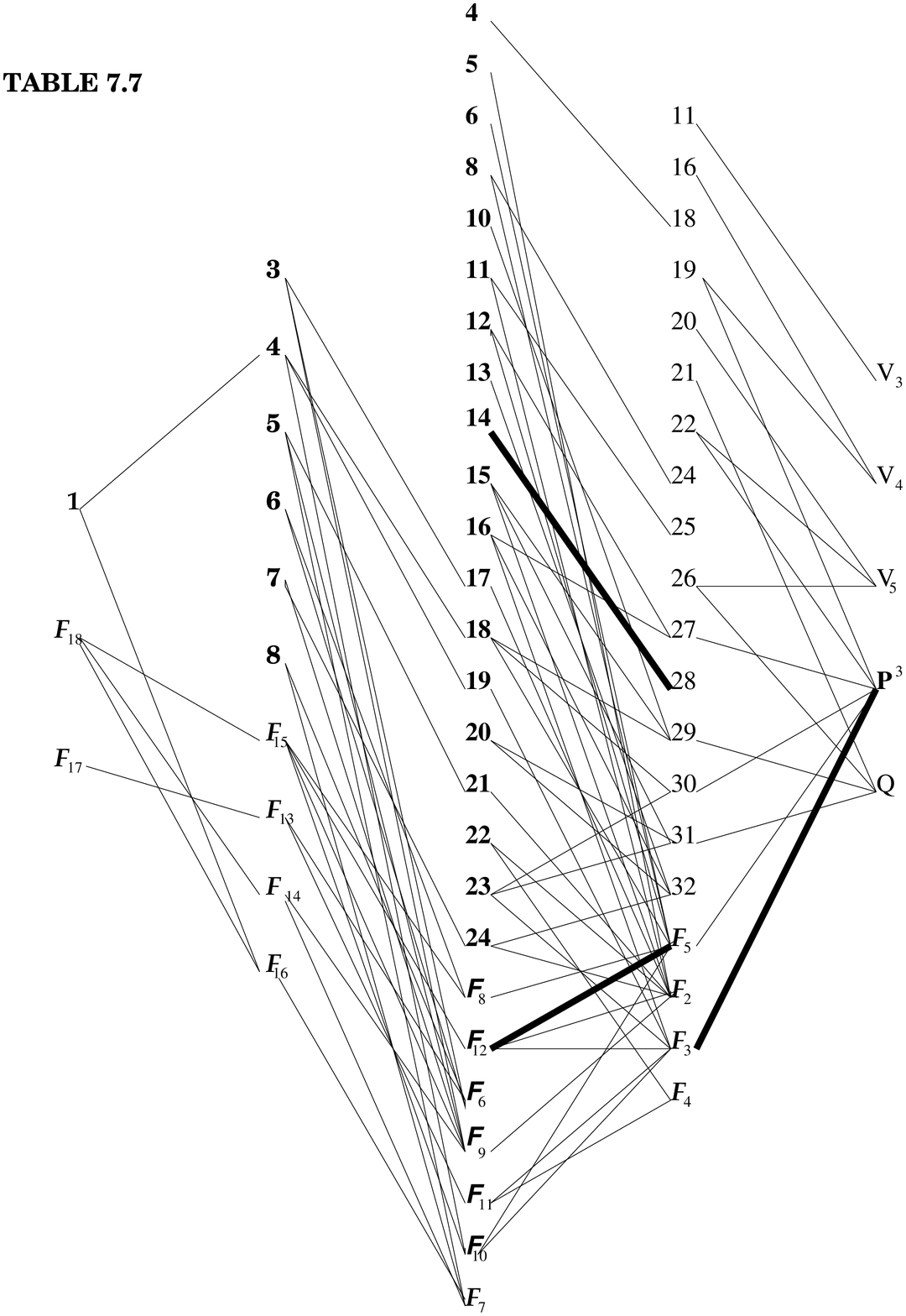}

\enddocument